\documentclass[11pt]{amsart}

\usepackage{amsfonts}
\usepackage{bm}
\usepackage{float}
\usepackage{geometry}
\usepackage{fancyhdr}
\usepackage{graphicx}

\allowdisplaybreaks[4]

\usepackage{tabularx}
\usepackage{bbm}
\usepackage{amsmath}
\usepackage{mathtools}
\usepackage{amsthm}
\usepackage{amssymb}
\usepackage{pdfpages}
\usepackage{listings}
\usepackage{setspace}
\usepackage{algorithm2e}
\usepackage{algpseudocode}
\DeclareMathOperator{\spanop}{span}
\newcommand{\rd}{\mathrm{d}}
\newcommand{\RR}{\mathbb{R}}
\newcommand{\NullL}{\mathrm{Null}\mathcal{L}}
\newcommand{\Vmat}{\mathsf{V}}
\newcommand{\Amat}{\mathsf{A}}
\newcommand{\Dmat}{\mathsf{D}}

\newtheorem{theorem}{Theorem}
\newtheorem{remark}{Remark}

\title[Coupling BGK model and Euler equations through linearized
Knudsen layer]{A numerical method for coupling the BGK model and Euler equations through the linearized Knudsen layer}

\author{Hongxu Chen}
\address{Mathematics Department, University of Wisconsin-Madison, 480 Lincoln Dr., Madison, WI 53705 USA.}
\email{hchen463@wisc.edu}

\author{Qin Li}
\address{Mathematics Department, University of Wisconsin-Madison, 480 Lincoln Dr., Madison, WI 53705 USA.}
\email{qinli@math.wisc.edu}

\author{Jianfeng Lu}
\address{Department of Mathematics, Department of
  Physics, and Department of Chemistry, Duke University, Box 90320,
  Durham NC 27708, USA}
\email{jianfeng@math.duke.edu}

\date{\today}

\thanks{The work of H.C.~and Q.L.~is supported in part by National Science Foundation under the grant
  DMS-1619778 and DMS-1750488. The work of J.L. is supported in part by the National
  Science Foundation under the grant DMS-1415939. The collaboration is
  also supported by National Science Foundation through KI-Net under
  grants RNMS-1107291 and RNMS-1107444. We thank Weiran Sun for
  helpful discussions. }

\begin{document}
\maketitle
\begin{abstract}
  The Bhatnagar-Gross-Krook (BGK) model, a simplification of the Boltzmann equation, in the absence of boundary effect, converges to
  the Euler equations when the Knudsen number is small. In practice, however, Knudsen layers emerge at the physical boundary, or at the
  interfaces between the two regimes. We model the Knudsen layer using a half-space kinetic equation, and apply a half-space numerical solver~\cite{lls:17_3, lls:15_2} to quantify the transition between the kinetic to the fluid regime. A full domain numerical solver is developed with a domain-decomposition approach, where we apply the Euler solver and kinetic solver on the appropriate subdomains and connect them via the half-space solver. In the nonlinear case, linearization is performed upon local Maxwellian. Despite the lack of analytical support, the numerical evidence nevertheless demonstrate that the linearization approach is promising.
\end{abstract}

\pagenumbering{arabic}

\section{Introduction}
The Bhatnagar-Gross-Krook (BGK) model, as a simplified version of the Boltzmann equation, is a classical model that describes the dynamics of rarified gas on the phase space. It has been extensively used in aerospace engineering, nuclear engineering and many related areas. In the dimensionless form the equation reads
\begin{equation}\label{eqn:BGK}
\partial_t F+v\cdot \nabla_xF=\frac{1}{\varepsilon}(M[F]-F),\quad  (t,x,v)\in \mathbb{R}^+\times\Omega\times \mathbb{R}^{d}\,,
\end{equation}
where $F=F(t,x,v)$ is the density function of the particles at time
$t\in \mathbb{R}^+$, position and velocity $(x,v)$ in the phase space
$\Omega\times\mathbb{R}^d$, where $\Omega\subset\mathbb{R}^d$ is the
physical domain and $\mathbb{R}^d$ is the velocity domain, and $d$ is
the spatial dimension. A dimensionless parameter $\varepsilon > 0$ is
called the Knudsen number representing the ratio of the mean free path
and the characteristic domain length scale. The term on the right hand
side of the equation is called the BGK term, with $M[F]$, termed the
Maxwellian distribution (local equilibrium), being a Gaussian with its
moments depending on $F$. More specifically:
\begin{equation}\label{eqn:Maxwellian}
  M[F]=\frac{\rho}{(2\pi T)^{d/2}}\exp\left(-\frac{|v-u|^2}{2T}\right)\,,
\end{equation}
with the macroscopic quantities $\rho$, $u$ and $T$ are obtained by taking the moments of $F$:
\begin{equation}\label{eqn:moments}
 \begin{pmatrix} \rho(x,t) \\ \rho(x,t) u(x,t) \\ d \rho(x,t) T(x,t)\end{pmatrix} =  \int_{\mathbb{R}^d} F(x,v,t)\begin{pmatrix}
1 \\v \\|v-u|^2\end{pmatrix} \rd{v}\,.
\end{equation}
Such definition enforces the Maxwellian to share the first $d+2$ moments with $F$, namely
\[
\int\phi(v)(M[F]-F)\rd{v} = 0\,,
\]
for $\phi(v) = 1,v$ or $|v|^2$ so that the entire equation conserves mass, momentum and energy.

We consider influx boundary conditions, meaning one is given the profile of $F$ on $\Gamma_-$, the collection of coordinates on the boundary of the physical domain with velocities pointing inside. Denoting $n_x$ the outer normal direction at $x\in\partial\Omega$, then
\[
F|_{\Gamma_-} = \phi(x,v)\,,\quad\text{with}\quad\Gamma_-=\{(x,v): x\in\partial\Omega\,, v\cdot n_x<0\}\,.
\]

The equation has drawn lots of attention on both theoretical and numerical aspects. On the analysis side, one can show that the equation, in the zero limit of $\varepsilon$, is asymptotically equivalent to the Euler equations, in the sense that $F$ converges to $M[F]$, and the macroscopic quantities follow the Euler equations. More specifically, the macroscopic quantities in~\eqref{eqn:Maxwellian} satisfy:
\begin{equation}\label{eqn:Euler}
\partial_t\mathcal{U} + \nabla_x\cdot\mathcal{F}(\mathcal{U}) = 0\,,
\end{equation}
with
\begin{equation*}
\mathcal{U}=\left(\rho\,,\rho u\,, E\right)\,,\quad\text{and}\quad \mathcal{F}(\mathcal{U})=\left(\rho u\,,\rho u\otimes u+\rho T\,, (E+\rho T)u\right)\,.
\end{equation*}
where $E=\frac{1}{2}\rho |u|^2+\frac{d}{2}\rho T$.

On the numerical side, the BGK equation, or in general, the
Boltzmann-like equation is very challenging. There are two main
numerical difficulties. First, the equation is on the phase
space instead of the physical space, so the discretization is done on
a higher dimensional space, leading to higher number of degree of freedoms. Second, the Knudsen number $\varepsilon$ could have many different
scales, and the solution behaves differently depending on which regime
the system is in. When $\varepsilon$ is relatively small, the
collision term is stiff, and to obtain numerical stability, the time
step has to resolve the small scale $\varepsilon$.

Extensive studies have been conducted to address the second
problem. In particular, ``asymptotic-preserving'' methods are
developed during the past decade. To a large extent, implicit
treatment of the equation needs to be employed to enlarge the stable
region, and thus relaxing the time step restrictions. However, for
kinetic equations, the stiff terms are typically nonlinear and
nonlocal, making implicit treatment complicated. A large body of work
is then proposed to find the surrogates of the collision term, or to
reformulate the equation to ``preserve'' the asymptotic
limit~\cite{Jin:99,FJ:10,LP:14,DP:11,Wild,BLM:08}, also see
reviews~\cite{DP:14,Jin_Rev:12}.

In this work, we are interested in tackling the first problem. In
particular, we aim at eliminating unnecessary dimensions whenever it
is possible. As argued above, asymptotically the BGK equation is
equivalent to the Euler equations in the fluid regime. Since the Euler
equations merely involve the spatial variables, in terms of reducing
computational complexity, one should compute the Euler equations,
instead of the BGK equation whenever it is a valid approximation. This approach is in
line with the domain decomposition method, as two different sets of
equations are applied in separate regions. However, there are some
immediate difficulties. Since the BGK equation and the Euler equations
are computed in different sub-domains and coupled at the interface, an
accurate coupling solver needs to be designed to correctly translate
data between the two systems.

The problem concerning the coupling of kinetic and fluid equations has
been a long-standing challenge. One major difficulty arises is the
boundary layer that emerges at the interface. At the interface,
$\varepsilon$ changes scales. Physically, it means the particles from
rarified regime suddenly are pushed into the condensed region, and it
takes a couple of mean free paths away from the sharp interface for
the particles to collide before the system achieving the local
equilibrium that is governed by the limiting fluid equation. This
drastic change is mathematically characterized by a boundary layer,
called the Knudsen layer in the current context, that characterizes the
damping from an arbitrary incoming rarified gas boundary condition to
a local equilibrium.

An intuitive coupling strategy is proposed in the pioneering work~\cite{LionsBP}, namely, one should compute the kinetic and fluid equations in separate domains, and the Knudsen layer equation as a model for the interface that translates data in between. Since the fluid solver and the kinetic solver are rather standard, the main challenge comes from the computation of the Knudsen layer equation. This approach was then used in some numerical studies, such as~\cite{Goudon2011, Golse1989, Klar1994, dellacherie2003coupling, dellacherie2003kinetic,golse1995numerical}, in most of which the layer equation was treated using Marshak condition~\cite{marshak1947note,Arnold1997}.

For the layer equation, there was a long stretch of investigation since 1970s. In the
linearized setting, firstly, it was shown in~\cite{bardos1986milne}
that the Knudsen layer equation is well-posed when the system has zero
bulk velocity (the so-called Milne problem), and then in the
celebrated paper~\cite{CGS:88} the authors extended the results and
were able to show the well-posedness for the linearized Boltzmann
equation with arbitrary bulk velocity when proper data is given, and
finally in~\cite{lls:17_3} the wellposedness was extended to a very
general class of linear kinetic systems under mild assumptions. This
particular work uses the weak formulation which makes the computation
feasible: with properly chosen basis functions, the Knudsen layer
equation becomes a coupled ODE system that gets rid of infinite domain
restriction. A spectral type algorithm was also proposed in the same
paper with rigorous error analysis that states quasi-optimality. We
emphasize that all these methods are made possible crucially depending
on the linearity. In the nonlinear regime, both the well-posedness and
numerical solver design are largely open, except a few results with
assumptions on weak nonlinearity and small data, as discussed
in~\cite{ukai2003nonlinear,ukai2004nonlinear}.

The mathematical challenge still persists today. Since the correct
boundary condition for a well-posed nonlinear Knudsen layer is still
unknown, no proper Knudsen layer solver has been proposed. It is not
our intention to address the nonlinearity in the Knudsen layer
equation. Rather, we would like to investigate, only numerically, if
linearization could in some sense capture the solution's nonlinear
behavior at the interface. We adopt the strategy proposed
in~\cite{LionsBP}, and at the interface, linearize the system upon a
suitably chosen Maxwellian function, hoping that in this small data
regime, the linearized Knudsen layer equation could serve as a good
approximation. This is a purely numerical approach, and we do not
intend to fully recover the nonlinear behavior. The aim is to propose
an intuitive and suitable strategy, and numerically investigate, to
what extent, the linearization is a good approximation.

Since the computation of the Knudsen layer is the main ingredient in the entire scheme, we first review it in Section 2. Section 3 and 4 are devoted to linearized and nonlinear setting of the coupling between the Knudsen layer computation and the AP solver used for regions without layers. Numerical evidence for linear and nonlinear setting will be demonstrated in the end of Section 3 and 4 respectively.

\section{Linearized systems and the Knudsen layer}\label{sec:set-up}
In this section we consider the linearization to the BGK equation and derive its asymptotic acoustic limit. We also give an overview of the theory and the numerical methods for the Knudsen layer equation in subsection~\ref{sec:layer}.
\subsection{Linearized BGK equation and its fluid limit}\label{sec:linearBE}
To perform linearization of the BGK equation~\eqref{eqn:BGK}, we assume the distribution function $F$ is close to a given global Maxwellian $M_\ast$ that has macroscopic state $(\rho_\ast,u_\ast,T_\ast)$, then substitute
\begin{equation}
  F=M_\ast+\sqrt{M_\ast}f,\quad\text{with}\quad M_\ast=\frac{\rho_\ast}{(2\pi T_\ast)^{d/2}}\exp\left(-\frac{|v-u_\ast|^2}{2T_\ast}\right)\,,
\end{equation}
in the BGK equation, with the first order expansion, we obtain the linearized BGK equation:
\begin{equation}\label{eqn:linear_BGK}
\begin{cases}
  \partial_t f+v\cdot \nabla_x f=\frac{1}{\varepsilon}\mathcal{L}_\ast f,\quad \text{with}\quad  \mathcal{L}_\ast f=m_\ast[f]-f\\
  f|_{\Gamma_-} = \frac{1}{\sqrt{M_\ast}}\left[\phi(x,v)-M_\ast\right]
  \end{cases}\,,
\end{equation}
where the linear Maxwellian $m_\ast[f](v)$ is a quadratic function that preserves the first $d+2$ moments of $f$ weighted by $\sqrt{M}_\ast$:
\begin{equation}\label{eqn:linear_moments_coll}
\langle f-m_\ast[f]\,,v^k\rangle_{M_\ast}= \int_\mathbb{R}(f-m_\ast[f])v^k\sqrt{M_\ast}\rd{v}=0,\quad k=0,1,2\,.
\end{equation}
Here $\langle \cdot,\cdot \rangle_{M_\ast}$ denotes the inner product with the weight $\sqrt{M_\ast}$:
\begin{equation}\label{eqn: weighted inner product}
  \langle f,g\rangle_{M_\ast}=\int_{\mathbb{R}}fg\sqrt{M_\ast}dv.
\end{equation}

For $d=1$, defining the moments $(\tilde{\rho},\tilde{u},\tilde{T})$ of $f$ by
\begin{equation}\label{eqn:macro_f}
\begin{pmatrix} \tilde{\rho} \\\tilde{\rho}u_\ast+\rho_\ast\tilde{u} \\\tilde{\rho}(u_\ast^2+T_\ast)+2\rho_\ast u_\ast\tilde{u}+\rho_\ast\end{pmatrix}=  \int_{\mathbb{R}}f
\begin{pmatrix}1 \\v \\v^2\end{pmatrix}\sqrt{M_\ast}\rd{v}=\langle f\,,
\begin{pmatrix}1 \\v \\v^2\end{pmatrix}\rangle_{M_\ast}\,,
\end{equation}
we can explicitly express $m_\ast(v)$:
\begin{equation}\label{eqn:linear_bgk_m}
  m_\ast[f](v)=\left[\frac{\tilde{\rho}}{\rho_\ast}+\frac{\tilde{u}}{T_\ast}(v-u_\ast)+\frac{\tilde{T}}{2T_\ast}\left(\frac{(v-u_\ast)^2}{T_\ast}-1\right)\right]\sqrt{M_\ast}\,.
\end{equation}

Formally, if $\varepsilon\to 0$ in~\eqref{eqn:linear_BGK}, the collision term $\mathcal{L}_\ast f$ dominates. It is a standard practice via the Hilbert expansion, without boundary layer effect, in the leading order, $f$ can be shown to be asymptotically equivalent to $m_\ast[f]$ whose macroscopic quantities satisfy the limiting acoustic equations:
\begin{equation}\label{eqn:acoustic}
\partial_t U+\Amat\cdot\partial_x U=0\,,
\quad \text{where}\quad
  \Amat=\begin{pmatrix}
        u_* & \rho_* & 0 \\
        \frac{T_*}{\rho_*} & u_* & 1 \\
        0 & 2T_* & u_* \\
      \end{pmatrix}
    \,,\quad \text{and}\quad U=\begin{pmatrix}\tilde{\rho} \\ \tilde{u} \\ \tilde{T}\end{pmatrix}\,.
\end{equation}

This limiting system is a hyperbolic system, and thus is diagonalizable with real eigenvalues. Writing $\Amat = \Vmat\cdot\Dmat\cdot\Vmat^{-1}$, we have:
\begin{equation}\label{eqn: Linearized Euleru}
\partial_t {U}+\Vmat\cdot \Dmat\cdot \Vmat^{-1}\cdot \partial_x {U}  = 0 \quad \Longrightarrow \quad \partial_t {\eta}+\Dmat\cdot \partial_x {\eta} = 0\,,
\end{equation}
where
\begin{equation}\label{eqn:acoustic_advection}
  \eta=\Vmat^{-1}\cdot U=\begin{pmatrix}                         \frac{T_\ast}{\rho_\ast}\tilde{\rho}-\frac{\tilde{T}}{2} \\
    \frac{\tilde{\rho}}{\rho_\ast}+\sqrt{\frac{3}{T_\ast}}\tilde{u}+\frac{\tilde{T}}{T_\ast} \\
    \frac{\tilde{\rho}}{\rho_\ast}-\sqrt{\frac{3}{T_\ast}}\tilde{u}+\frac{\tilde{T}}{T_\ast} \end{pmatrix}\,,\quad  \text{and}\quad   \Dmat=\textrm{diag}\Big(u_\ast,u_\ast+\sqrt{3T_\ast},u_*-\sqrt{3T_\ast}\Big)\,.
\end{equation}
Therefore $\eta_i$ satisfies the advection equation with speed $d_i=\Dmat_{ii}$ (for $i=1,2,3$). We note that quantities $\eta_i$ can also be directly obtained by taking the moments of $f$:
\begin{equation}\label{eqn:eta_projection}
\eta_i = \langle f\,, p_i\rangle_{M_\ast}\,,
\end{equation}
where
\begin{align*}
p_1 =& -\frac{1}{2\rho_\ast}(v-u_\ast)^2 + \frac{3T_\ast}{2\rho_\ast}\,, \\
p_{2,3} =& \frac{1}{\rho_\ast T_\ast}(v-u_\ast)^2 \pm \frac{1}{\rho_\ast}\sqrt{\frac{3}{T_\ast}}(v-u_\ast)\,.
\end{align*}

\subsection{Half-space kinetic boundary layer equation}\label{sec:layer}

The Knudsen layer equation was initially proposed in an early work~\cite{LionsBP} to describe the behavior of Knudsen layer. Suppose one is given a kinetic equation in a bounded domain with small Knudsen number, a boundary layer would emerge to translate the incoming boundary condition to a local equilibrium function whose macroscopic quantities are governed by the limiting fluid equations. The equation for this layer, termed the Knudsen layer equation, is formed by locally ``stretching" the coordinate and balance the leading order terms:
\begin{equation}\label{eqn:layer}
  \begin{cases}
     v\partial_z f = \mathcal{L}_\ast f, &  (z,v)\in\RR^+\times\RR\\
    f|_{z=0}=\phi(v;x,t), & v>0
  \end{cases}\,.
\end{equation}
In the equation, $z$ is the rescaled spatial variable in the layer with $z=0$ being mapped to the boundary coordinate $x\in\partial\Omega$, and $z=\infty$ being mapped to the interior of the domain along the negative ray of $n_x$ from $x\in\partial\Omega$. We note that the equation is a steady state problem with time $t$ and the boundary coordinate $x$ serving as parameters via the boundary condition, and $\phi$ comes from confining $f|_{\Gamma_-}$ to $x$.

Interested readers are referred to~\cite{LionsBP} for the derivation and we omit the details from here. In the rest of this section we summarize the well-posedness result and a spectral type algorithm for the Knudsen layer equation~\eqref{eqn:layer}.

\subsubsection{Wellposedness}
The Knudsen layer equation has a unique solution in $L_2$ only if certain boundary condition is satisfied at the limit, seen from the following theorem.
\begin{theorem}[Theorem 1 from~\cite{CGS:88} and Theorem 3 from~\cite{lls:17_3}]\label{thm:wellpose}
  Let the incoming data $\phi\in L^2((1+|v|)\mathbf{1}_{v>0}\rd{v})$,
  the half-space equation
  \begin{equation}
    \begin{cases}
      v\partial_z f = \mathcal{L}_\ast f, &  (z,v)\in\RR^+\times\RR\\
      f|_{z=0}=\phi(v), & v>0
    \end{cases}
  \end{equation}
  has a unique solution such that
  \begin{equation}
    \lim_{z \to \infty} f(z, \cdot) \in H^+\oplus H^0\,.
  \end{equation}
  Here $H^+$ and $H^0$ are the collections of positive and zero modes
  associated with multiplicative operator $v$ in $\NullL_\ast$, the
  null space of the collision operation $\mathcal{L}_\ast$.
\end{theorem}

The theorem is proved for general kinetic equations, and the definitions of $H^{+,0}$ are rather vague. In the following, we derive the explicit expression for the two spaces. According to the definition of the linearized BGK operator given in~\eqref{eqn:linear_bgk_m}, setting $\mathcal{L}_\ast[f] = 0$ naturally makes $f$ a quadratic function (weighted by $\sqrt{M_\ast}$), meaning:
\begin{equation}
  \NullL_\ast =\spanop \bigl\{\sqrt{M_\ast}\,,v\sqrt{M_\ast}\,,v^2\sqrt{M_\ast}\bigr\}\,.
\end{equation}
We rearrange this space using the following basis functions:
\begin{equation}\label{eqn:basis_NullL}
  \begin{cases} \displaystyle
    \chi_0(v)=\frac{1}{\sqrt{6\rho_\ast}}
    \Big(\frac{(v-u_\ast)^2}{T_\ast}-3\Big)\sqrt{M_\ast}  \\
    \displaystyle \chi_+(v)=\frac{1}{\sqrt{6\rho_\ast}}\Big(\sqrt{\frac{3}{T_\ast}}
    (v-u_\ast)+\frac{(v-u_\ast)^2}{T_\ast}\Big)\sqrt{M_\ast} \\
    \displaystyle \chi_-(v)=\frac{1}{\sqrt{6\rho_\ast}}
    \Big(\sqrt{\frac{3}{T_\ast}}(v-u_\ast)-\frac{(v-u_\ast)^2}{T_\ast}\Big)\sqrt{M_\ast}
  \end{cases}\,.
\end{equation}

It is a good set of basis functions due to the following three properties it satisfy:
\begin{itemize}
\item[1.] they form an orthogonal expansion of the null space, namely
  \[
 \NullL_\ast = \spanop\{\chi_0\,,\chi_+\,,\chi_-\}\,,\quad\text{and}\quad \int\chi_i(v)\chi_j(v)\rd{v} =\langle\chi_i\,,\chi_j\rangle = \delta_{ij}\,,\quad \forall i,j \in \{0, +, -\}\,;
  \]
\item[2.] they are eigenfunctions of the multiplicative operator $v$
  restricted in $\NullL_\ast$:
  \begin{equation*}
    \langle v\chi_i\,,\chi_j\rangle = u_i\delta_{ij}\,;
  \end{equation*}
\item[3.] the eigenvalues $u_i$ are given by the bulk velocity and the temperature:
  \begin{equation}\label{eqn:speed_modes}
    u_0=u_*,\quad u_+=u_*+\sqrt{3T_*},\quad u_-=u_*-\sqrt{3T_*}\,.
  \end{equation}
\end{itemize}
Here $\langle \cdot,\cdot\rangle$ represents the standard inner product:
\begin{equation}\label{eqn: inner product}
\langle f,g\rangle=\int_{\mathbb{R}}fgdv.
\end{equation}
Since that $H^+$ and $H^0$ are eigensubspace of the multiplicative operator $v$ restricted on $\NullL_\ast$, one has:
\begin{equation*}
H^+=\spanop\{\chi_i|u_i>0\},\quad H^-=\spanop\{\chi_i |u_i<0\},\quad\text{and}\quad H^0=\spanop\{\chi_i|u_i=0\}\,,
\end{equation*}
For convenience of the notation, define the dimension of the spaces $\nu_{\pm,0} = \text{dim}H^{\pm,0}$ and re-label the modes, we have:
\begin{equation}\label{eqn:def_zeta}
H^\pm=\spanop\{\zeta_{\pm,1}\,,\cdots,\zeta_{\pm,\nu_\pm}\},\quad H^0=\spanop\{\zeta_{0,1}\,,\cdots,\zeta_{0,\nu_0}\}\,.
\end{equation}
According to Theorem~\ref{thm:wellpose}, for uniqueness, at $z=\infty$, the solution has to be in:
\begin{equation*}
    \lim_{z \to \infty} f(z, \cdot) \in  \spanop\{\zeta_{0,1}\,,\cdots,\zeta_{0,\nu_0}\,,\zeta_{+,1}\,, \cdots,\zeta_{+,\nu_+}\} = \spanop\{\chi_i|u_i\geq 0\}\,.
\end{equation*}
\begin{remark}
We note that the theorem discusses the uniqueness and asserts that the projection of $f|_{z=\infty}$ on $H^-$ should be zero. If we remove this requirement, the equation loses its uniqueness but we still have the existence. In fact, one can show there are infinitely many solutions. The solution space is simply:
\begin{equation*}
f + H^-
\end{equation*}
where $f$ is the unique solution in Theorem~\ref{thm:wellpose}.
\end{remark}

\subsubsection{Spectral method for half-space equations}
Finding the numerical solution to the half-space problem, however, carries a different challenging aspect. The problem is supported on an infinite domain, and cannot be numerically treated easily. In~\cite{lls:17_3}, a semi-analytic spectral method was developed that achieves quasi-optimality in $v$ and is analytic in $x$. The algorithm relies on the damping-recovering approach. To be more specific, a damping term is introduced and added to the right hand side of the equation, so that all elements in $\NullL_\ast$ are damped out from the solution, forcing the solution to the damped equation, denoted by $f_d$, to be zero at $z=\infty$. The trick of finding the solution to the original equation lies in the fact that a very special boundary condition can be designed, so that when it is put into the damped equation, it cancels the effect of the damping term.

More explicitly, the damped Knudsen layer equation reads:
\begin{equation}\label{eqn:damped}
\begin{cases}
v\partial_z f_d=\mathcal{L}_\ast f_d + \mathcal{L}_df_d;\\
f_d(z=0,v)=\phi(v)\,,\quad v>0; \\
\lim_{z \to \infty} f_d(z, \cdot ) = 0.
\end{cases}
\end{equation}
where the added damping term is:
\begin{align*}
 \mathcal{L}_df=\sum_{i\in\{+,-,0\}}\sum^{\nu_i}_{j=1}\alpha (v+u_\ast)\zeta_{i,j}\langle v\zeta_{i,j}\,,f\rangle+\sum_{i=1}^{\nu_0}\alpha (v+u_\ast)\mathcal{L}^{-1}((v+u_\ast)\zeta_{0,j})\langle v\mathcal{L}^{-1}((v+u_\ast)\zeta_{0,j})\,,f\rangle\,,
\end{align*}
with some small arbitrarily chosen value for $\alpha$ and its value does not affect the results. We summarize the recovering formula in the following theorem.

\begin{theorem}[Proposition 3.4 from~\cite{lls:17_3}]\label{thm:boundary_compute}
Let $\phi\in L^2((1+|v|)\mathbf{1}_{v>0}\rd{v})$, then
\begin{equation}\label{eqn:soln}
f=f_d-\sum_{i=1}^{\nu_+}\xi_{+,i}(g_{+,i}-\zeta_{+,i})-\sum_{j=1}^{\nu_0}\xi_{0,j}(g_{0,j}-\zeta_{0,j})
\end{equation}
uniquely solves~\eqref{eqn:layer} with $\lim_{z \to \infty} f=f_{\infty}\in H^+\oplus H^0$ being the end-state:
\begin{equation}\label{eqn:end_state}
f_{\infty}=\sum_{j=1}^{\nu_+}\xi_{+,j}\zeta_{+,j}+\sum_{k=1}^{\nu_{0}}\xi_{0,k}\zeta_{0,k}=\sum_{j=1}^{\nu_+}\langle f_\infty\,,\zeta_{+,j}\rangle\zeta_{+,j}+\sum_{k=1}^{\nu_{0}}\langle f_\infty\,,\zeta_{0,k}\rangle\zeta_{0,k}\,.
\end{equation}
Here $f_d$ solves~\eqref{eqn:damped} with inflow data $\phi$, and $g_{i,j}$ solves~\eqref{eqn:damped} with inflow data $\zeta_{i,j}$. $\{\xi_{+,j},\xi_{0,k}\}$ are coefficients that solve:
\begin{equation*}
C\cdot \vec{\xi} = \vec{Q}\,,\quad\text{with}\quad C=\left(\begin{array}{cc}C_{++}\,,C_{+0}\\C_{0+}\,,C_{00}\end{array}\right)\,,\quad\vec{\xi} = \left(\begin{array}{c}\vec{\xi}_{+}\\ \vec{\xi}_{0}\end{array}\right)\,,\quad\text{and}\quad\vec{Q} = \left(\begin{array}{c}\vec{Q}_{+}\\ \vec{Q}_{0}\end{array}\right)\,,
\end{equation*}
where the vector coefficients $\vec{\xi}_{+}=[\xi_{+,1}\,,\xi_{+,2}\,,\cdots,\xi_{+,\nu_+}]^T$, and $\vec{Q}_{+}=[\langle v\zeta_{+,1}\,,f_d\rangle\,,\cdots\langle v\zeta_{+,\nu_+}\,,f_d\rangle]^T$ (and similarly for $\vec{\xi}_0$ and $\vec{Q}_0$), and the matrix is defined by:
\begin{align*}
C_{++,ij} = \langle v\zeta_{+,i}\,,g_{+,j}|_{z=0}\rangle\,,\quad C_{+0,ij} = \langle v\zeta_{+,i}\,,g_{0,j}|_{z=0}\rangle\,,\\
C_{0+,ij} = \langle v\zeta_{0,i}\,,g_{+,j}|_{z=0}\rangle\,,\quad C_{00,ij} = \langle v\zeta_{0,i}\,,g_{0,j}|_{z=0}\rangle\,.
\end{align*}
\end{theorem}

The theorem suggests the following steps to compute $f$:
\begin{itemize}
\item[1.] compute the damped equation~\eqref{eqn:damped} for $f_d$ and $g_{i,j}$ using the associated boundary conditions;
\item[2.] use $g_{i,j}$ to define the matrices $C$ and use $f_d$ to define $\vec{Q}$ for computing $\vec{\xi}$;
\item[3.] assemble $f$ according to~\eqref{eqn:soln} and $f|_{z=\infty}$ according to~\eqref{eqn:end_state}.
\end{itemize}
The details of the algorithm that computes the damped equation~\eqref{eqn:damped} are summarized in Appendix~\ref{sec:app_A}.

\section{Acoustic limit of the linearized BGK equation with kinetic boundary condition}\label{sec:linear}
With a spectral accurate numerical solver for the Knudsen layer equation, we are now ready to numerically couple the layer equation and the interior fluid equation. We investigate the coupling in the pure linearized regime in this section and leave the nonlinear coupling to Section 4.

To demonstrate the coupling, we consider one particular example where the linearized BGK system with small Knudsen number governing a finite bounded domain with a non-Maxwellian type incoming flow. More specifically:

\begin{equation}\label{eqn:linear_system2}
    \begin{cases}
      \partial_t f+v\partial_x f=\frac{1}{\varepsilon}\mathcal{L}_\ast f = \frac{1}{\varepsilon}(m_\ast-f), &   (t,x,v)\in \mathbb{R}^+\times[0,1]\times \mathbb{R}\,,\quad\varepsilon\ll 1\\
      f(t=0,x,v)=f_0(x,v)\\
      f(t,x=0,v)=\phi_l(t,v), & v>0\\
      f(t,x=1,v)=\phi_r(v,t), & v<0
    \end{cases}\,.
\end{equation}
As $\varepsilon\to 0$, for any interior domain $[a,b]\subsetneqq[0,1]$, this equation is well approximated by the linearized Euler equations and its diagonalization. We rewrite~\eqref{eqn: Linearized Euleru} into
\begin{equation}
\partial_t U+\Amat\cdot\partial_x U=0\,,\quad\Rightarrow\quad \partial_t\eta_i +d_i\partial_x\eta_i=0\,.
\end{equation}
Depending on the sign of $d_i$, $\eta_i$ is either advecting left or right.

The two sets of equations share the following properties:
\begin{itemize}
\item The speed $d_i$ are counterparts of the ``averaged speed'' of $\zeta_{\pm/0}$ defined in~\eqref{eqn:def_zeta}, namely (with an arbitrary ordering)
\begin{equation*}
d_1 = u_0 = u_\ast\,,\quad d_2 = u_+ = u_\ast+\sqrt{3T_\ast}\,,\quad d_3 = u_- = u_\ast-\sqrt{3T_\ast}\,.
\end{equation*}
\item The projection of $f\vert_{z=\infty}$ on $\zeta$, denoted as $\xi$ as in Theorem~\ref{thm:boundary_compute}, is a counterpart of $\eta$. If we compare the definition of $\eta$ in~\eqref{eqn:eta_projection} and $\chi_{0,\pm}$ in~\eqref{eqn:basis_NullL}, we see that:
\begin{equation*}
\chi_0 = -\frac{2\sqrt{\rho_\ast M_\ast}}{\sqrt{6}T_\ast}p_1\,,\quad\chi_+ = \frac{\sqrt{\rho_\ast M_\ast}}{\sqrt{6}}p_2\,,\quad\chi_- = -\frac{\sqrt{\rho_\ast M_\ast}}{\sqrt{6}}p_3\,.
\end{equation*}
and thus $\zeta$s, the re-labels of $\chi$s are $p$s multiplied by constants, meaning $\xi = \langle f\,,\zeta\rangle$ and $\eta = \langle f\,,p\rangle_{M_\ast}$ are counterparts of each other, with the specific matching determined by the Mach number.
\end{itemize}

These transformations are crucial in translating the inflow boundary condition on $f$ to the Dirichlet type boundary condition on $d_i$. Indeed, with $d_i>0$, $\eta_i$ would be a right-propagating mode with its Dirichlet boundary condition imposed on $x=0$, and that piece of information should come from $\phi_l$. Meanwhile, the Knudsen layer equation supported close to $x=0$ would be projecting information in $H^+$, whose dimension exactly depends on the signs of $\{u_{0},u_+,u_-\}$.

\subsection{Numerical method}
The numerical method is straightforward, and we discuss it briefly. For a simpler presentation, below we assume the Mach number is between $0$ and $1$ so that $u_\ast+\sqrt{3T_\ast}>u_\ast>0$ and $u_\ast-\sqrt{3T_\ast}<0$. In this case, in the fluid regime, $\eta_1$ and $\eta_2$ are right-propagating modes needing information from $x=0$, while $\eta_3$ is the left-propagating mode needing information from $x=1$. The method works similarly for other cases of Mach number with straightforward adjustment.

We evenly discretize the domain into $N$ cells with cell length $h = \frac{1}{N}$:
\begin{equation*}
0 = x_0 < x_1<x_2<\cdots<x_N=1\,,
\end{equation*}
and we denote $\eta^n_{i,j}$ the numerical estimate of $\eta_i(t_n,x_j)$. With upwinding scheme, one has:
\begin{align}\label{eqn:update_eta}
&\eta_{i,j}^{n+1}=\eta_{i,j}^n+\frac{u_i\Delta t}{h}\left(\eta_{i,j}^n-\eta_{i,j-1}^n\right)\,,\quad &i=1,2,\quad j=2,\cdots,N \,;\\
&\eta_{3,j}^{n+1}=\eta_{3,j}^n+ \frac{u_3\Delta t}{h}\left(\eta_{3,j+1}^n-\eta_{3,j}^n\right)\,,\quad &j=1,\cdots,N-1\,.\nonumber
\end{align}
This updates all $\eta^n_{ij}$ except the three boundary points: $\eta_{1,0}$, $\eta_{2,0}$ and $\eta_{3,N}$.

To update the boundary condition for $\eta$, one needs to compute the Knudsen layer equation with kinetic boundary inflow $\phi_l$ and $\phi_r$. Take the left boundary for example, since $u_\ast>0>u_\ast-\sqrt{3T_\ast}$, at $x=0$,
\[
H^+=\spanop\{\chi_0\,,\chi_+\}\quad\text{and}\quad H^-=\spanop\{\chi_-\}\,.
\]
The information restricted on $H^-$ should be provided by the left-going mode $\eta_{3}(x=0)$, while $\eta_{1/2}(x=0)$ need to be computed according to the layer equation. To be more specific, realizing
\begin{align*}
\begin{cases}
\xi_{-,1} = \langle f^l|_{z=\infty}\,,\zeta_{-,1}\rangle = \langle f^l|_{z=\infty}\,,\chi_-\rangle = \frac{-\sqrt{\rho_\ast}}{\sqrt{6}}\langle f|_{z=\infty}\,,p_3\rangle_{M_\ast} = \frac{-\sqrt{\rho_\ast}}{\sqrt{6}}\eta_3(x=0)\,,\\
\xi_{+,1} = \langle f^l|_{z=\infty}\,,\zeta_{+,1}\rangle = \langle f^l|_{z=\infty}\,,\chi_0\rangle = \frac{-2\sqrt{\rho_\ast}}{\sqrt{6}T_\ast}\langle f|_{z=\infty}\,,p_1\rangle_{M_\ast} = \frac{-2\sqrt{\rho_\ast}}{\sqrt{6}T_\ast}\eta_1(x=0)\,,\\
\xi_{+,2} = \langle f^l|_{z=\infty}\,,\zeta_{+,2}\rangle = \langle f^l|_{z=\infty}\,,\chi_+\rangle = \frac{\sqrt{\rho_\ast}}{\sqrt{6}}\langle f^l|_{z=\infty}\,,p_2\rangle_{M_\ast} = \frac{\sqrt{\rho_\ast}}{\sqrt{6}}\eta_2(x=0)\,,
\end{cases}
\end{align*}
we first subtract the $H^-$ mode provided by $\eta_3(x=0)$ from the inflow data:
\[
\phi^l_\text{modify} = \phi_l(v) - \frac{\sqrt{\rho_\ast}}{\sqrt{6}}\eta^{n+1}_{3,0}\chi_{-}(v)\,,\quad v>0\,,
\]
and the modified Knudsen layer equation reads:
\begin{equation}\label{eqn:layer2}
\begin{cases}
     v\partial_z f^l = \mathcal{L}_\ast f^l\,,\quad (z,v)\in\RR^+\times\RR\,,\\
    f^l|_{z=x=0}=\phi^l_\text{modify}(v)\,.
\end{cases}
\end{equation}
The boundary layer equation \eqref{eqn:layer2} is then solved using the procedure in Theorem~\ref{thm:boundary_compute} for $f^l_{z=\infty}$, and one updates the boundary condition for $\eta$:
\begin{equation}\label{eqn:eta_bdry_l}
 \eta^{n+1}_{1,0} = \frac{-\sqrt{6}T_\ast}{2\sqrt{\rho_\ast}}\langle f^l|_{z=\infty},\chi_0\rangle  \,,\quad  \eta^{n+1}_{2,0} = \frac{\sqrt{6}}{\sqrt{\rho_\ast}}\langle f^l|_{z=\infty},\chi_+\rangle \,.
\end{equation}

The same procedure can be done for the right boundary at $x=1$ with
$z=\frac{1-x}{\varepsilon}$ and $\tilde{v} = -v$. Since the sign is
flipped, we have $H^+ = \spanop\{\chi_-\}$ and $H^- = \spanop\{\chi_0\,,\chi_+\}$, and solver the follow layer equation:
\begin{equation}\label{eqn:layer3}
\begin{cases}
     \tilde{v}\partial_z f^r = \mathcal{L}_\ast f^r\,,\quad (z,v)\in\RR^+\times\RR\,,\\
    f^r|_{z=0 \,(x=1)}=\phi_r(-v) + \frac{\sqrt{\rho_\ast}}{\sqrt{6}}\eta^{n+1}_{1,N}\chi_{+}(v)  - \frac{2\sqrt{\rho_\ast}}{\sqrt{6}T_\ast}\eta^{n+1}_{2,N}\chi_{0}(v)\,,\quad v>0\,.
\end{cases}
\end{equation}
Then $\eta_3$ at the right end is updated accordingly:
\begin{equation}\label{eqn:eta_bdry_r}
\eta^{n+1}_{3,N} = \frac{-\sqrt{6}}{\sqrt{\rho_\ast}}\langle f^r|_{z=\infty},\chi_-\rangle\,.
\end{equation}

We summarize the algorithm below.

\RestyleAlgo{boxruled}
\begin{algorithm}[H]
\SetAlgoLined
\KwData{\begin{itemize}\item[1] Kinetic boundary condition $\phi_l(v,t_{n+1})$ for $v>0$ and $\phi_r(v,t_{n+1})$ for $v<0$;\item[2] Solution at $t_n$: all $\eta^n_{i,j}$ for $i=1,2,3$, and $j =0,1,\cdots,N$.
\end{itemize}}
\KwResult{Solution at $t_{n+1}$: all $\eta^{n+1}_{i,j}$}
Step I: Update $\eta^{n+1}_{i,j}$ in the interior according to~\eqref{eqn:update_eta};\\
Step II: Compute~\eqref{eqn:layer2} and~\eqref{eqn:layer3} using the procedure from Theorem~\ref{thm:boundary_compute} for $f^{l/r}|_{z=\infty}$;\\
Step III: Update boundary data:\\
     \For{$i=1,2,3$}{
     \eIf{$u_i >0$}{
     Set $\eta_{i,0}$ according to~\eqref{eqn:eta_bdry_l};
     }{Set $\eta_{i,N}$ according to~\eqref{eqn:eta_bdry_r};}
     }
     \caption{Updating the computation of the acoustic limit from $t_n$ to $t_{n+1}$}
     \end{algorithm}
\subsection{Numerical examples}
We demonstrate several numerical examples for computing the acoustic limit for the linearized BGK equation. In all the tests we have, domain is set to be $[0,1]$. To compute the BGK equation, Knudsen number is chosen to be $\varepsilon=\{\frac{1}{32},\frac{1}{64},\frac{1}{128},\frac{1}{256}\}$. Numerically, velocity domain is chosen to be $v\in[-16,16]$ large enough so that the Gaussian tales excluded is negligible. In space, $h = 10^{-3}$ so that the layers are resolved, and $\Delta t = h/20$ to satisfy the CFL condition. To compute the acoustic limit we follow the Algorithm 1, and use sample grids using $h = 0.005$ and $\Delta t = h/5$. To measure the error, we set
\begin{equation}\label{eqn:error_numerics}
D_{{\rho}} = \|\tilde{\rho}_\text{BGK} - \tilde{\rho}_\text{Euler}\|_{L_2([0.1,0.9])}
\end{equation}
with ``BGK" referring to the solution of the BGK equation, and ``Euler" indicating the solution to the linearized Euler equations with Knudsen layer correction. $D_{{u}}$ and $D_{{T}}$ are computed similarly.

\begin{figure}[htb]
  \centering
  \includegraphics[width=19cm, height=5cm]{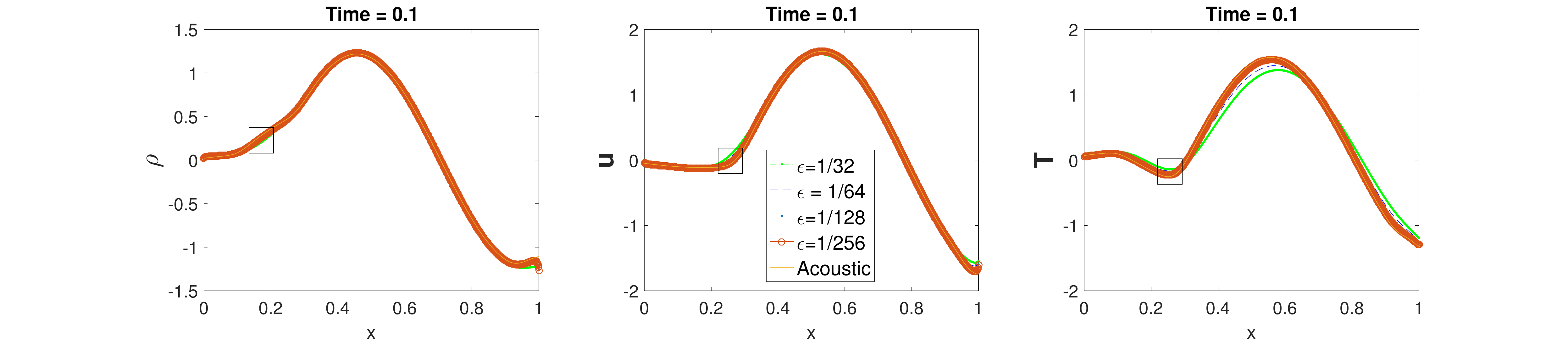}\\
  \caption{Test 1: Solution at Time$=0.1$.}
  \label{fig:Test1}
\end{figure}

\begin{figure}[htb]
  \centering
  \includegraphics[width=19cm, height=5cm]{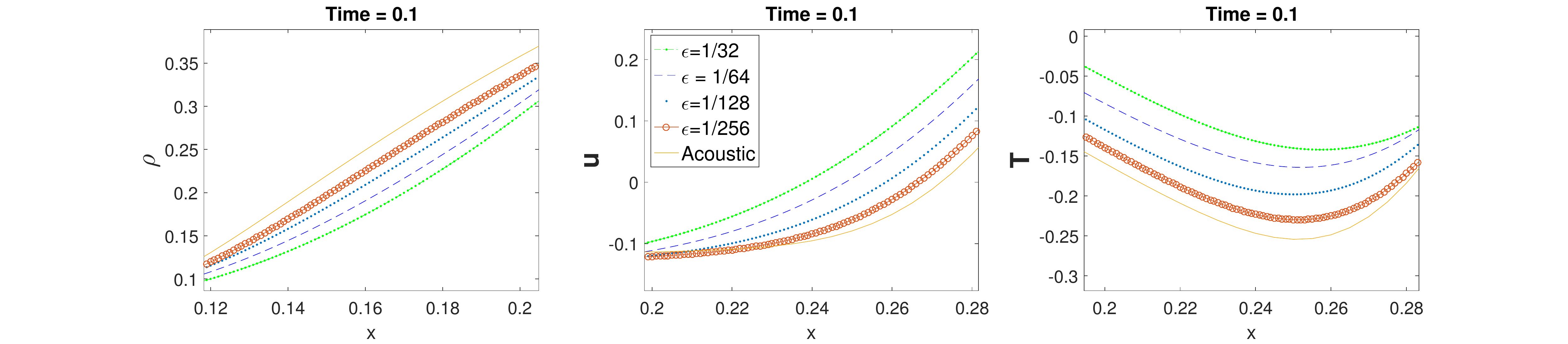}\\
  \caption{Test 1: Zoomed in the box region in Figure~\ref{fig:Test1}.}
  \label{fig:Test1zoom}
\end{figure}

\begin{figure}[htb]
  \centering
  \includegraphics[width=19cm, height=5cm]{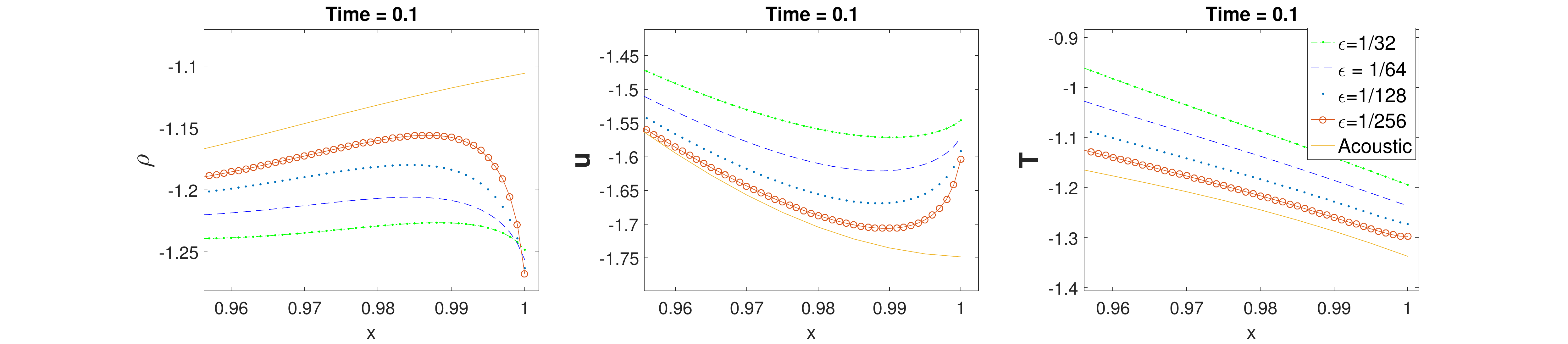}\\
  \caption{Test 1: Boundary layer zoomed in at $x=1$.}
  \label{fig:Test1zoom2}
\end{figure}

\begin{figure}[htb]
  \centering
  \includegraphics[width=19cm, height=5cm]{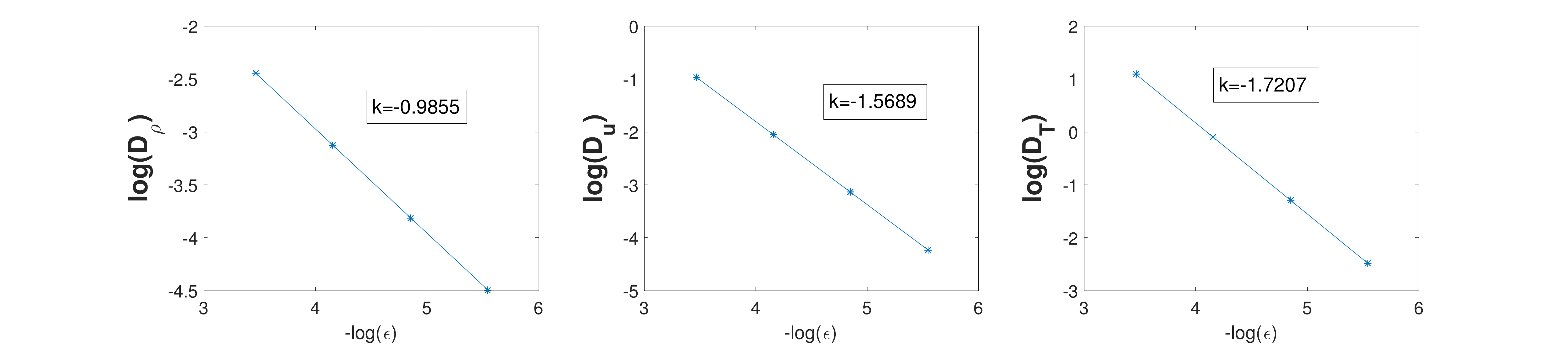}\\
  \caption{Test 1: Error $D_{{\rho}}$, $D_{{u}}$ and $D_{{T}}$ as functions of $\varepsilon$ in log-log scale.}
  \label{fig:Test1error}
\end{figure}

\textbf{Test 1: Subsonic case with boundary layers emerging at both ends.} In this test, we set $\rho_\ast=1,u_\ast=1,T_\ast=1$, and thus the problem is the subsonic regime with $u_\ast+c_\ast>u_\ast>0>u_\ast-c_\ast$ with $c_\ast = \sqrt{3T_\ast}$. The initial condition and boundary condition are given by:
\begin{equation*}
  \left\{
    \begin{array}{ll}
      \text{Boundaries: } f_l(v,t)=0, \quad  f_r(v,t)=0  \\
      \text{Initial: } f_0(x,v)=\frac{1.5\sin(2\pi x)}{\rho_\ast}+\frac{1.5\sin(2\pi x)}{T_\ast}(v-u_\ast)+\frac{1.5\sin(2\pi x)}{2T_\ast}(\frac{(v-u_\ast)^2}{T_\ast}-1)
    \end{array}
  \right.\,. 
\end{equation*}
Correspondingly, we can derive the initial and boundary condition for the acoustic limit. More specifically, using~\eqref{eqn:linear_bgk_m}, we have initial data:
\[\tilde{\rho}_0(x,v)=1.5\sin(2\pi x),\quad  \tilde{u}_0(x,v)=1.5\sin(2\pi x),\quad \tilde{T}_0(x,v)=1.5\sin(2\pi x).\]
To compute the boundary condition for the acoustic limit, we apply Algorithm 1 using~\eqref{eqn:eta_bdry_l} with $\phi_l=0$. For this particular case, denote $c_-=(c_{-,0},c_{-,+})$, the solution of the left layer solution using $\chi_-$ as the incoming data, projected onto the corresponding modes, namely let:
\begin{equation*}
v\partial_z f+\mathcal{L}_\ast f=0\,,\quad f|_{z=0}=\chi_-(v)\,,\quad v>0,
\end{equation*}
then we set $c_-=(c_{-,0},c_{-,+})$ the coefficients for $f|_{z=\infty}=c_{-,0}\chi_0+c_{-,+}\chi_+$. Similarly, solve the right layer solution using $\chi_{+,0}$ as the incoming data:
\begin{equation*}
v\partial_zf +\mathcal{L}_\ast f=0\,, \quad f|_{z=0}=\chi_{+/0}(-v)\,,\quad (v>0),
\end{equation*}
and denote $c_{+/0}$ are coefficients: $f|_{z=\infty}=c_{+/0}\chi_-$. Then the boundary condition for acoustic equations at $x=0$ and $x=1$ can be made explicit:
\[
\eta^{n+1}_{1,0}=\frac{\eta^{n+1}_{3,0}}{2}c_{-,0}\,,\quad \eta^{n+1}_{2,0}=-\eta^{n+1}_{3,0}c_{-,+}\,,\quad \eta^{n+1}_{3,N}=2\eta^{n+1}_{1,N}c_0-\eta^{n+1}_{2,N}c_+\,.\]
In Figure~\ref{fig:Test1}, we plot the solution at $\text{Time}=0.1$. The distribution function starts with a Maxwellian and to this point, different profiles with different $\epsilon$ have not deviated from each other too much. However, with a zoomed-in profile, discrepancy is shown, see Figure~\ref{fig:Test1zoom}. The layer profile is significant, as demonstrated in Figure~\ref{fig:Test1zoom2}, in which it is clear that smaller $\epsilon$ leads to better approximation to the acoustic limit. We document the error and show its dependence on $\epsilon$ on log-log scale in~\ref{fig:Test1error}: the decay is a clear straight line. This is a numerical evidence that the error decays algebraically fast in the linear setting. Particularly, the decay of $D_\rho$ is about $\mathcal{O}(\epsilon)$.

\textbf{Test 2: Supersonic case, with boundary layer emerging at $x=1$ only. Compatible initial and boundary data.} In the second test, we set the reference state as $\rho_\ast=1,u_\ast=2,T_\ast=1/2$. The left boundary conditions are placed in $\NullL_\ast$ to avoid the left boundary layer. The initial data and boundary data are set to be compatible. This is a supersonic case and all boundary conditions for the acoustic limit should be placed at the left end. The initial and boundary conditions for the linearized BGK equation are given:
\begin{equation*}
  \left\{
    \begin{array}{ll}
      \text{Boundaries: } f_l(v,t)=t\chi_++t\chi_0+t\chi_-, \quad  f_r(v,t)=0  \\
      \text{Initial: } f_0(x,v)=\frac{1.25\sin(2\pi x)}{\rho_\ast}+\frac{1.25\sin(2\pi x)}{T_\ast}(v-u_\ast)+\frac{1.25\sin(2\pi x)}{2T_\ast}(\frac{(v-u_\ast)^2}{T_\ast}-1)
    \end{array}
  \right.\,. 
\end{equation*}
We note the initial and boundary conditions are compatible in the sense that
\[f_l(v,0)=f_0(0,v)=0,\quad f_r(v,0)=f_0(1,v)=0.\]
Correspondingly, we can compute for the initial and boundary conditions for the acoustic limit:
\begin{equation*}
  \left\{
    \begin{array}{ll}
      \text{Boundaries}:
      \eta (x=0)=\left(-\frac{\sqrt{6}T_\ast}{2\sqrt{\rho_\ast}}t\,,\frac{\sqrt{6}}{\sqrt{\rho_\ast}}t\,,-\frac{\sqrt{6}}{\sqrt{\rho_\ast}}t\right)^\top\\
      \text{Initial: }\tilde{\rho}_0(x,v)=1.25\sin(2\pi x),\quad  \tilde{u}_0(x,v)=1.25\sin(2\pi x),\quad \tilde{T}_0(x,v)=1.25\sin(2\pi x)
    \end{array}
  \right.\,. 
\end{equation*}

\begin{figure}[htb]
  \centering
  \includegraphics[width=19cm, height=5cm]{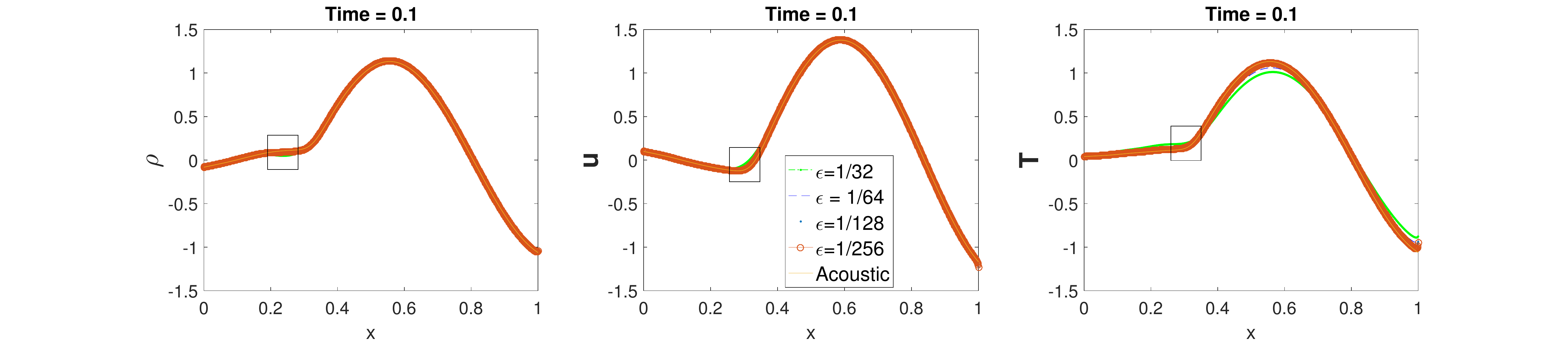}\\
  \caption{Test 2: Solution at Time$=0.1$.}
  \label{fig:Test2}
\end{figure}

\begin{figure}[htb]
  \centering
  \includegraphics[width=19cm, height=5cm]{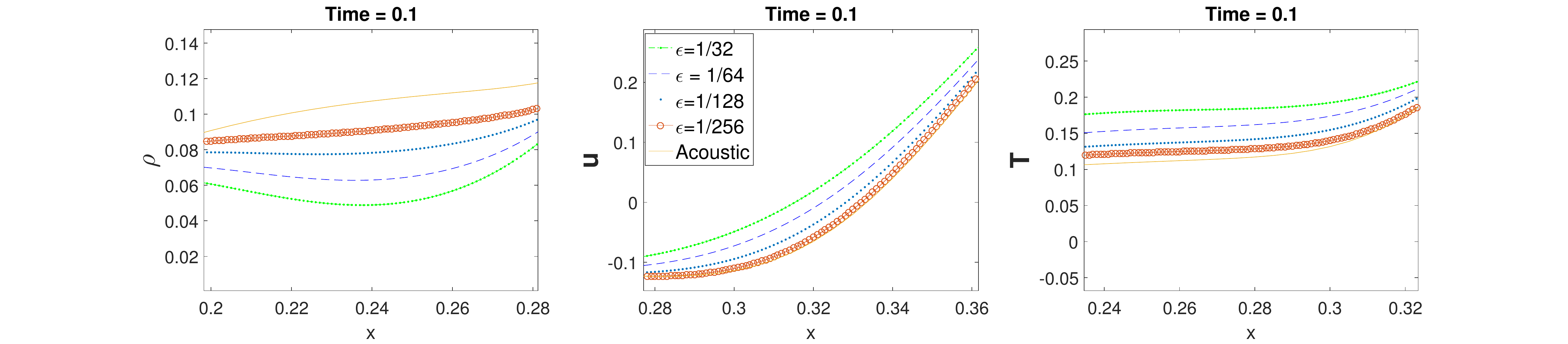}\\
  \caption{Test 2: Zoomed in the box of Figure~\ref{fig:Test2} to see the approximation.}
  \label{fig:Test2zoom}
\end{figure}

\begin{figure}[htb]
  \centering
  \includegraphics[width=19cm, height=5cm]{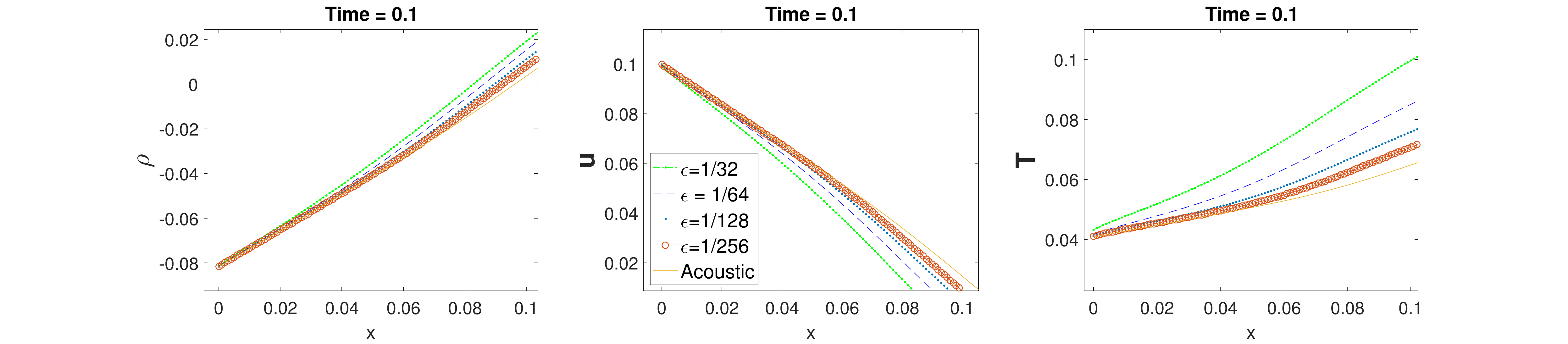}\\
  \caption{Test 2: Left boundary zoomed at $x=0$.}
  \label{fig:Test2zoom2}
\end{figure}

\begin{figure}[htb]
  \centering
  \includegraphics[width=19cm, height=5cm]{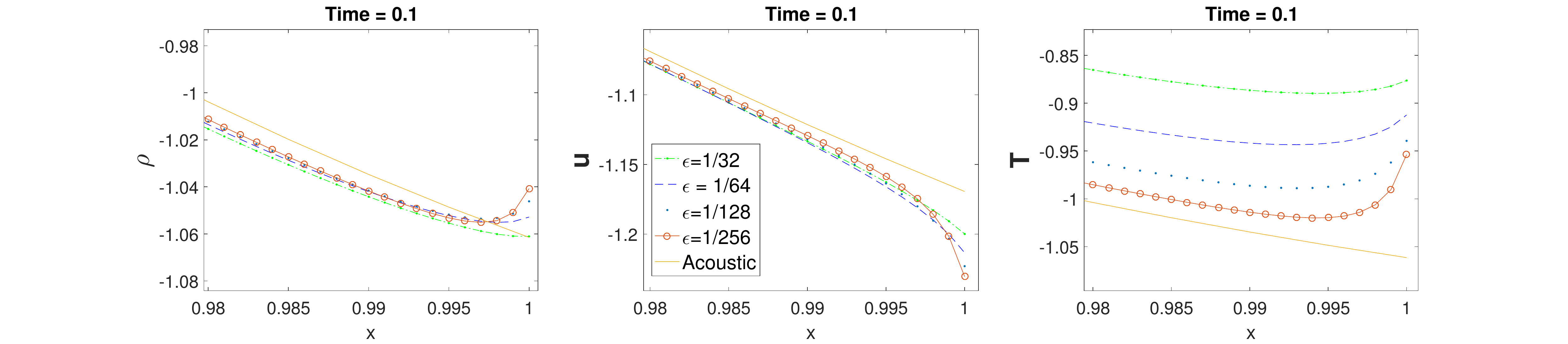}\\
  \caption{Test 2: Boundary layer zoomed at $x=1$.}
  \label{fig:Test2zoom3}
\end{figure}

\begin{figure}[htb]
  \centering
  \includegraphics[width=19cm, height=5cm]{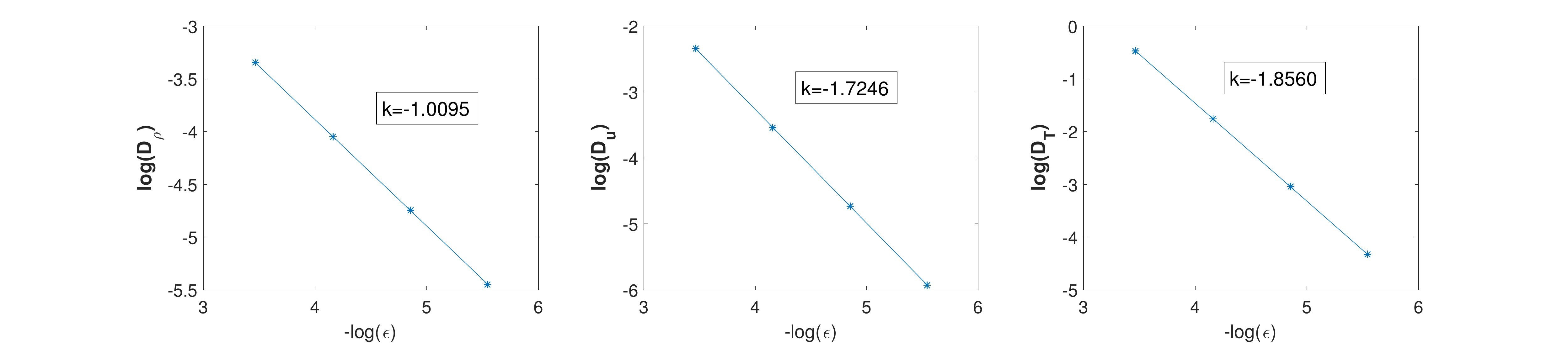}\\
  \caption{Test 2: Error $D_{{\rho}}$, $D_{{u}}$ and $D_{{T}}$ as functions of $\varepsilon$ in log-log scale.}
  \label{fig:Test2error}
\end{figure}

In Figure~\ref{fig:Test2} and Figure~\ref{fig:Test2zoom}, we plot the solution and a small region zoomed-in at $\text{Time}=0.1$. No layer occurs at $x=0$ as shown in Figure~\ref{fig:Test2zoom2}, due to the carefully chosen boundary condition. One Knudsen layer does emerge at $x=1$, and the profile is plotted in Figure~\ref{fig:Test2zoom3}. Smaller $\epsilon$ leads to closer approximation to the acoustic limit. We also document the errors, and plot them in log-log scale with respect to $\epsilon$, as seen in Figure~\ref{fig:Test2error}.

\begin{figure}[htb]
  \centering
  \includegraphics[width=19cm, height=5cm]{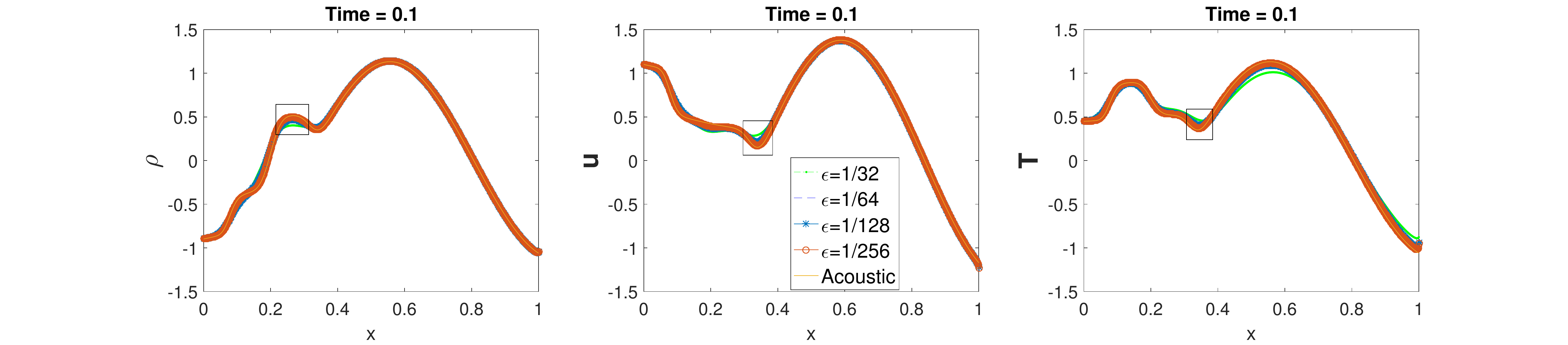}\\
  \caption{Test 3: Solution at Time$=0.1$.}
  \label{fig:Test3}
\end{figure}

\begin{figure}[htb]
  \centering
  \includegraphics[width=19cm, height=5cm]{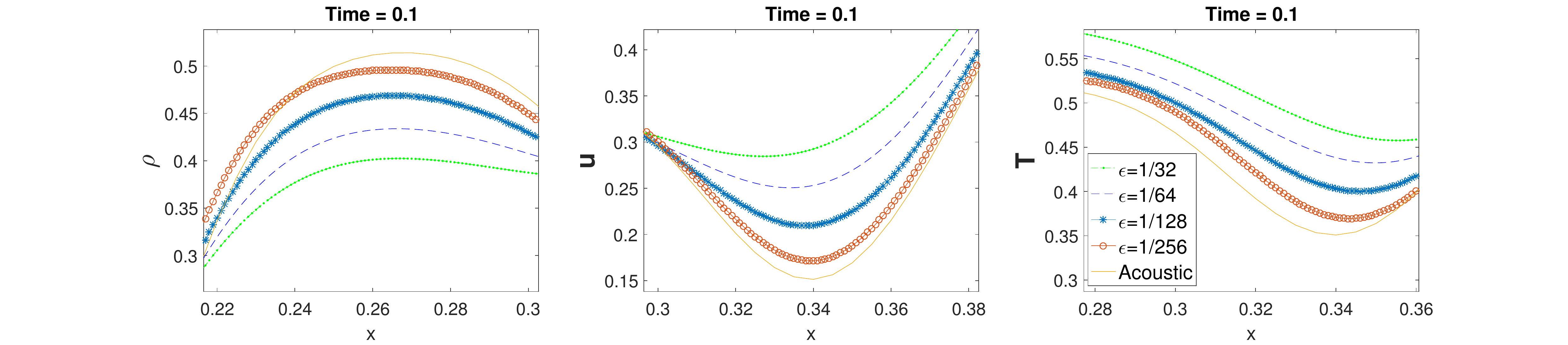}\\
  \caption{Test 3: Zoomed in the box of Figure~\ref{fig:Test3} to see the approximation.}
  \label{fig:Test3zoom}
\end{figure}

\begin{figure}[htb]
  \centering
  \includegraphics[width=19cm, height=5cm]{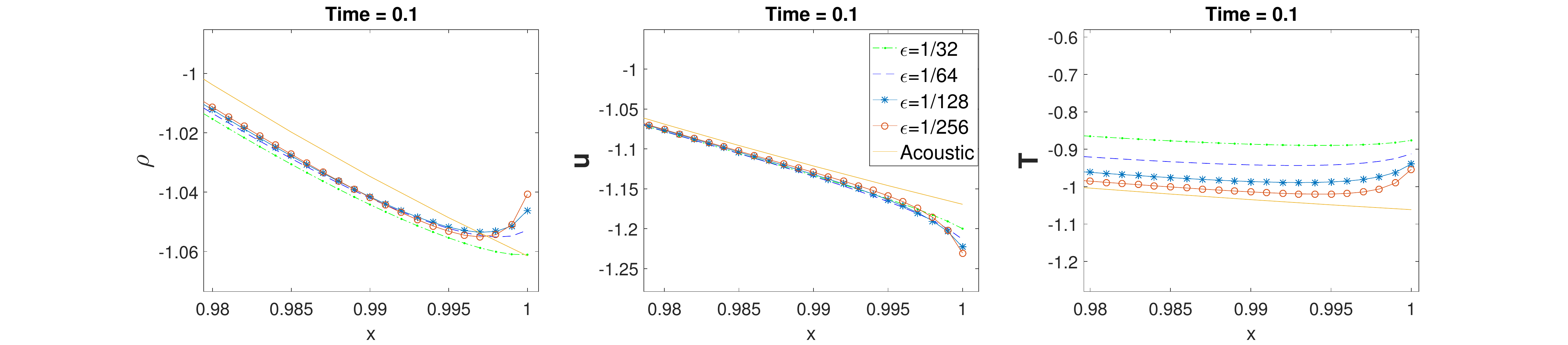}\\
  \caption{Test 3: Boundary layer zoomed at $x=1$.}
  \label{fig:Test3zoom2}
\end{figure}

\begin{figure}[htb]
  \centering
  \includegraphics[width=19cm, height=5cm]{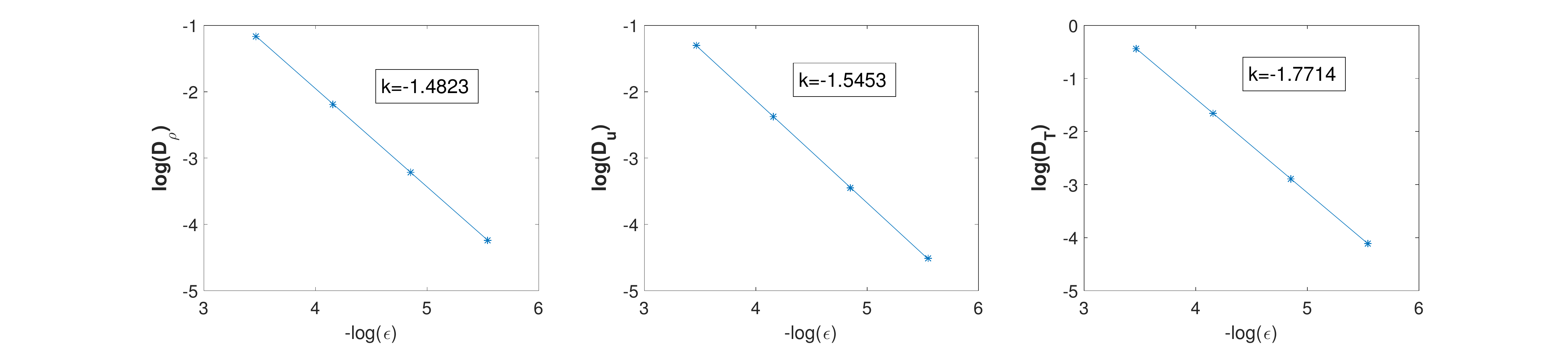}\\
  \caption{Test 3: Error $D_{{\rho}}$, $D_{{u}}$ and $D_{{T}}$ as functions of $\varepsilon$ in log-log scale.}
  \label{fig:Test3error}
\end{figure}

\textbf{Test 3: Supersonic case, with boundary layer emerging at $x=1$. Initial condition is incompatible with boundary data.} In the third test, we use the same parameters by setting $\rho_\ast=1,u_\ast=2,T_\ast=1/2$. However we adjust the initial data and boundary data to be incompatible. The initial and boundary condition given to the BGK equation are:
\begin{equation}\label{}
  \left\{
    \begin{array}{ll}
      \text{Boundaries: } f_l(v,t)=(1+t)\chi_++(1+t)\chi_0+(1+t)\chi_-,\quad f_r(v,t)=0 \\
      \text{Initial: }  f_0(x,v)=\frac{1.25\sin(2\pi x)}{\rho_\ast}+\frac{1.25\sin(2\pi x)}{T_\ast}(v-u_\ast)+\frac{1.25\sin(2\pi x)}{2T_\ast}(\frac{(v-u_\ast)^2}{T_\ast}-1)
    \end{array}
  \right. .
\end{equation}
The initial and boundary conditions are incompatible in the sense that
\[f_l(v,0)\neq f_0(0,v),\quad f_r(v,0)\neq f_0(1,v).\]
Correspondingly we have the data for the acoustic limit:
\begin{equation*}
  \left\{
    \begin{array}{ll}
      \text{Boundaries}:\eta(x=1) =\left(-\frac{\sqrt{6}T_\ast}{2\sqrt{\rho_\ast}}(1+t)\,,\frac{\sqrt{6}}{\sqrt{\rho_\ast}}(1+t)\,,-\frac{\sqrt{6}}{\sqrt{\rho_\ast}}(1+t)\right)^\top\\
 \text{Initial: }\tilde{\rho}_0(x,v)=1.25\sin(2\pi x),\quad  \tilde{u}_0(x,v)=1.25\sin(2\pi x),\quad \tilde{T}_0(x,v)=1.25\sin(2\pi x)
    \end{array}
  \right. .
\end{equation*}
Similar as the previous two examples, in Figure~\ref{fig:Test3} and Figure~\ref{fig:Test3zoom}, we plot the solution and a small region zoomed-in at $\text{Time}=0.1$. The Knudsen layer emerges at $x=1$ but not at $x=0$, and the layer is plotted in Figure~\ref{fig:Test3zoom2}. The errors' decay with respect to $\epsilon$ is plotted in Figure~\ref{fig:Test3error}. Still smaller $\epsilon$ leads to closer approximation to the acoustic limit. We emphasize that comparing of Figure~\ref{fig:Test3error} with Figure~\ref{fig:Test2error}, we clearly see that the error is significantly larger: in Test 3 log-scale $D_\rho$ ranges from $-1$ to $-4$ while that in Test 2 ranges from $-3$ to $-6$.

\section{Coupling the BGK and the Euler equations}\label{sec:nonlinear}
In this section, we numerically investigate the real challenging problem: we
study the coupling between the nonlinear BGK equation and the
nonlinear Euler equations. More specifically, we perform linearization
at the boundaries at every time step to obtain a linearized Knudsen
layer equation, and compute the flux to exchange information between
kinetic and fluid solvers. We largely follow the strategy proposed
in~\cite{Goudon2011}, however, unlike making extra assumptions to the layer equation as was done in~\cite{Goudon2011}, here we do compute the half-space Knudsen layer equation with the half-space solver that ensures spectral accuracy. The entire process has two main sources of error: 1. higher
order terms are thrown away due to the linearization procedure,
introducing the linearization error; 2. numerical error is introduced
during the computation of the half-space Knudsen layer equation. The main difference between the current studies and the previous results in~\cite{dellacherie2003coupling, dellacherie2003kinetic,golse1995numerical, Klar1994} is that we would like to completely eliminate the second type of error. For that we apply the recently developed
half-space Knudsen layer solver. It has the spectral accuracy in
velocity domain and is analytic in spatial variable. With this error reduced, the error in our
computation essentially only comes from the linearization, allowing us to truly see to what extent could linearization approximate the nonlinear Knudsen layer.

We emphasize that linearization is, at the current stage, the only solution to such nonlinear kinetic equations with layers, largely due to the lack of well-posedness theory on the analytical level. We do not claim the proposed algorithm below is the optimal choice, but rather, numerically test to what extent can linearization approximate the nonlinear coupling.

We would be focusing on the following set-up:
\begin{equation}\label{eqn:BGK_nonlinear}
\begin{cases}
&\partial_t F + v\partial_xF = \frac{1}{\varepsilon}(M[F] - F)\,,\quad (t,x,v) \in\mathbb{R}^+\times[0,1]\times\mathbb{R}\\
&F(t,x=0,v) = F_l(t,v)\,,\quad v>0\\
&F(t,x=1,v) = F_r(t,v)\,,\quad v<0
\end{cases}
\end{equation}
where
\[
\varepsilon \begin{cases}
= 1\,,\quad x\in[1/2,1]\\
\ll 1\,,\quad x\in[0,1/2]
\end{cases}\,.
\]

Since the system is in fluid regime in subdomain $[0,1/2]$, we expect two boundary layers emerging at the two ends of this subdomain, namely, the layer would appear at $x=0$, the physical boundary, and $x=1/2$, the interface.

The nonlinear nature makes this problem significantly harder. Since there is no ``global-Maxwellian", ``local-Maxwellian" needs to be found at each time step, upon which linearization is performed for us to obtain the linearized Knudsen layer equation.

For a numerical setup, we divide the domain into $N$ cells and apply Finite Volume type method. We denote the grids $x_{k}$ the cell centers, and fluxes are then computed on the half-grids $x_{k+1/2}$:
\begin{equation}
0 = x_{1/2} < x_{3/2}<\cdots<x_{N/2+1/2} = 1/2 < \cdots <x_{N-1/2}<x_{N+1/2}= 1\,.
\end{equation}
In the velocity domain we truncate the computational domain to be $[-16,16]$ and use $32M$ evenly distributed grid points for the velocity discretization.
\begin{equation}\label{eqn:vj}
-16 = v_{0} < v_{1}<\cdots<v_{16M} = 0 < \cdots <v_{32M-1}<v_{32M}= 16\,.
\end{equation}
The velocity cut-off is chosen to be big enough for mass to be almost
conserved in our simulations.

The Euler equations will be computed in $x<1/2$ and the BGK equation will be computed in $x>1/2$ with a to-be-specified Knudsen layer equation computed in the middle $x=1/2$ to couple them. The numerical methods for both the Euler equations and the BGK equation are rather standard, and we briefly review them in Section 4.1. The computation of the Knudsen layer equation will be discussed in Section 4.2. We summarize the algorithms before demonstrating numerical examples in the end of this section.

\subsection{Numerical methods for the Euler and the BGK equations}
Finite volume method will be applied to treat the limiting compressible Euler equations~~\eqref{eqn:Euler} in numerical domain $[0,1/2]$.

Denote $\mathcal{U}_i^n$ the numerical solution to the equation in cell $[x_{i-1/2},x_{i+1/2}]$ at time $t_n$, and $\mathcal{F}^n_{i+1/2}$ the numerical flux at the cell boundary $x=x_{i+1/2}$. From time step $t_n$ to $t_{n+1}$, one has:
\begin{equation}\label{eqn:Roe method}
\mathcal{U}_i^{n+1}=\mathcal{U}_i^n -\frac{\Delta t}{h}(\mathcal{F}^n_{i+1/2}-\mathcal{F}^n_{i-1/2})\,,
\end{equation}
with fluxes prepared at each cell border $\mathcal{F}^n_{i-1/2}$ for all $i=1,\cdots, N/2+1$.

For the flux term in the interior with $i=2\,,\cdots,N$, we follow the standard Roe flux method and choose
\begin{equation}\label{eqn:flux_interior}
\mathcal{F}_{i-1/2}=\frac{1}{2}[\mathcal{F}(\mathcal{U}_{i-1})+\mathcal{F}(\mathcal{U}_i)]-\frac{1}{2}|\hat{A}_{i-1/2}|(\mathcal{U}_i-\mathcal{U}_{i-1})\,,\quad i = 2,3,\cdots N/2\,,
\end{equation}
where $\hat{A}_{1-1/2}$ is the Jacobian of $\nabla_\mathcal{U}\mathcal{F}$:
\begin{equation*}
\hat{A}_{i-1/2}=\left(
                    \begin{array}{ccc}
                      0 & 1 & 0 \\
                      \frac{1}{2}(\gamma-3)\hat{u}^2 & (3-\gamma)\hat{u} & \gamma-1 \\
                      \frac{1}{2}(\gamma-1)\hat{u}^3-\hat{u}\hat{H} & \hat{H}-(\gamma-1)\hat{u}^2 & \gamma\hat{u} \\
                    \end{array}
                  \right)
\end{equation*}
evaluated using averaged velocity $\hat{u}$, total specific enthalpy $\hat{H}$ and sound speed $\hat{c}$:
\begin{align*}
\begin{cases}
\hat{u}&=\frac{\sqrt{\rho_{i-1}}u_{i-1}+\sqrt{\rho_i}u_i}{\sqrt{\rho_{i-1}}+\sqrt{\rho_i}}\\
\hat{H}&=\left[(E_{i-1}+\rho_{i-1}T_{i-1})/\sqrt{\rho_{i-1}}+(E_i+\rho_iT_i)/\sqrt{\rho_i}\right]/\left[\sqrt{\rho_{i-1}}+\sqrt{\rho_i}\right]\\
\hat{c}&=\sqrt{2(\hat{H}-\frac{1}{2}\hat{u}^2)}
\end{cases}\,.
\end{align*}

For the flux terms at the two ends of the domain, $\mathcal{F}_{1/2}$ and $\mathcal{F}_{N/2+1/2}=0$ are still unknown. Formula~\eqref{eqn:flux_interior} stops being valid at the boundaries and boundary condition needs to be incorporated. We defer the discussion to Section~\ref{sec:flux}.

To compute the BGK equation, we adopt the approach taken in~\cite{coron1991numerical}. With the pre-set discretization, we denote $F^n_{i,j}$ the numerical approximation at $(x_i,v_j)$ at time step $t^n$. The computation is split into two steps:
\begin{align*}
\partial_tF + v\partial_xF = 0\,,\quad \partial_tF = M[F]-F\,,
\end{align*}
with the correspondingly schemes:
\begin{eqnarray}\label{eqn:update_kinetics}
F^{n+1/2}_{i,j} &= F^n_{i,j}  - \Delta tv_j(\partial_xF)^n_{i,j}\,,\quad F^{n+1}_{i,j} &= F^{n+1/2}_{i,j} - {\Delta t}(M^{n+1/2}_{i,j}-F^{n+1/2}_{i,j})\,.
\end{eqnarray}
where the Maxwellian in the second part of the scheme is defined as:
\begin{equation*}
  M^{n+1/2}_{i,j}=\frac{\rho^{n+1/2}_{i,j}}{(2\pi T^{n+1/2}_{i,j})^{1/2}}\exp\left(-\frac{(v-u^{n+1/2}_{i,j})^2}{2T^{n+1/2}_{i,j}}\right)\,,
\end{equation*}
with its moments computed by:
\begin{equation*}
\left(\rho^{n+1/2}\,, \rho^{n+1/2} u^{n+1/2}\,,E^{n+1/2}\right)^\top = \sum_j\Delta vF_{ij}^{n+1/2}\left(1\,,v_j\,,v_j^2/2\right)^\top\,.
\end{equation*}

To compute the transport part in~\eqref{eqn:update_kinetics}, the simple upwinding method is used, namely:
\begin{align}\label{eqn:Nonlinear BGK step1}
\begin{cases}
  F_{i,j}^{n+1/2}=F_{i,j}^n-\frac{v_j\Delta t}{h}\left(F_{i,j}^n-F_{i-1,j}^n\right)\,, & \quad\hbox{for $i=N/2+2,\cdots,N+1$ and $v_j>0$}\,; \nonumber\\
F_{i,j}^{n+1/2}=F_{i,j}^n- \frac{v_j\Delta t}{h}\left(F_{i+1,j}^n-F_{i,j}^n\right)\,, & \quad\hbox{for $i=N/2+1,\cdots,N$ and $v_j<0$}\,.
\end{cases}
\end{align}
In the formulation $F^n_{N+1,j}$ takes the boundary condition $F_r(t_n,v_j)$ with $v_j<0$. This process leaves out the update for $F^{n+1/2}_{N/2+1,j}$ for $v_j>0$, and it would require information at the interface from the Euler equations' side. The Knudsen layer equation is computed, as will be discussed below.

\subsection{Boundary flux at the interface}\label{sec:flux}
From the analysis above, it is seen clearly that there are three terms that need the boundary information:

\begin{itemize}
\item $\mathcal{F}_{1/2}$ and $\mathcal{F}_{N/2+1/2}$, the fluxes at the interfaces for the Euler equations;
\item $F_{N/2+1,j}$ for $v_j>0$, the incoming flow for the kinetic equation.
\end{itemize}

These terms will be determined by finding a good Maxwellian function to perform linearization upon, and computing the linearized Knudsen layer equation. Below we discuss the computation of $\mathcal{F}_{1/2}$ as an example. The computation of the other two terms is similar.

According to the definition of the flux:
\begin{equation}\label{eqn:flux}
 \mathcal{F}(\mathcal{U})= \left(\rho u\,,\rho u^2+\rho T \,,(E+\rho T)u\right)=  \int v\left(1\,,v\,,|v|^2/2\right)F\rd{v}\,,
\end{equation}
so to numerically obtain $\mathcal{F}_{1/2}$, one needs to find $F_{1/2,j}$ for all $j$. Assume at time step $t_n$, the distribution function $F$ around $x_{1/2} = 0$ is close to the local Maxwellian, denoted by $M_\ast$. Defining the reference macroscopic state $(\rho_\ast\,, u_\ast\,,T_\ast)$, and writing
\begin{equation}\label{eqn:nonlinearF_expansion}
F = M_\ast + \sqrt{M_\ast}f\,,
\end{equation}
one derives that $f$ satisfies the Knudsen layer equation (ignoring higher order terms and stretching coordinates $z=\frac{x}{\varepsilon}$):
\begin{equation}\label{eqn:f_layer_nonlinear}
\begin{cases}
&v\partial_zf = m_\ast - f\,,\quad (t,z,v) \in\mathbb{R}^+\times[0,\infty]\times\mathbb{R}\\
&f(z=0,v) = \frac{F-M_\ast}{\sqrt{M_\ast}}(t,x=0,v) = \frac{F_l - M_\ast}{\sqrt{M_\ast}}\,,\quad v>0\\
\end{cases}\,,
\end{equation}
with $m_\ast$ is the linear infinitesimal determined by $(\rho_\ast\,, u_\ast\,,T_\ast)$. As analyzed in the previous section, $z=\infty$ corresponds to the end of the layer, which can be regarded as $x=x_{1/2}=0$. At this point, $f\in\NullL$, meaning:
\begin{equation}\label{eqn:f_infty}
f = f_- + f_+\,,\quad\text{with}\quad f_- = \sum_{i=1}^{\nu_-}\xi_{-,i}\zeta_{-,i}\,,\quad f_+ = \sum_{i=1}^{\nu_0}\xi_{0,i}\zeta_{0,i}+\sum_{i=1}^{\nu_+}\xi_{+,i}\zeta_{+,i}\,
\end{equation}
in charge of the flow sending inwards from the wall to the interior. Plugging~\eqref{eqn:nonlinearF_expansion} and~\eqref{eqn:f_infty} into~\eqref{eqn:flux}, one gets:
\begin{equation}\label{eqn:boundary_flux}
  \mathcal{F}_{1/2}=\int \left[M_\ast+\sqrt{M_\ast}f_-+\sqrt{M_\ast}f_+\right]v\left(\begin{array}{c}1 \\v \\|v|^2/2\end{array}\right)\rd{v}\,,
\end{equation}
where $f_+$ is in charge of the information getting into the interior from the physical boundary that includes the positive modes in the Knudsen layer equation, and $f_-$ is in charge of the information flowing out of the domain. According to Theorem~\ref{thm:boundary_compute}, the computation of $f_+$ is standard once $f_-$ is known.

To find $f_-$ that is consistent with the local Maxwellian function selected, we firstly define the fluctuation in macroscopic quantities:
\begin{equation}\label{eqn:fluctuation}
U_{\textrm{fluc}}=U_1-U_\ast=(\rho_{\textrm{fluc}},u_{\textrm{fluc}},T_{\textrm{fluc}})^\top\,,
\end{equation}
and its associated infinitesimal Maxwellian:
\begin{equation}\label{eqn:fluctuation inf}
m_{U_{\textrm{fluc}}}=\left[\frac{\rho_{\textrm{fluc}}}{\rho_\ast}+u_{\textrm{fluc}}\frac{v-u_\ast}{T_\ast}+\frac{T_{\textrm{fluc}}}{2T_\ast}\left(\frac{(v-u_\ast)^2}{T_\ast}-1\right)\right]\sqrt{M_\ast}\,.
\end{equation}
We then determine $f_-$ by projecting this local infinitesimal onto the negative modes:
\begin{equation}\label{eqn:f_-}
f_-(v)=\sum_{i=1}^{\nu_-}\xi_{-,i} \zeta_{-,i}\,,\quad\text{with}\quad \xi_{-,i}=\frac{\int vm_{U_\textrm{fluc}} \chi_{-,i} \rd{v}}{\int v|\chi_{-,i}|^2 \rd{v}}\,.
\end{equation}

Running the algorithm presented in Theorem~\ref{thm:boundary_compute} for $f_+$ and utilize formula~\eqref{eqn:boundary_flux}, one obtains $\mathcal{F}_{1/2}$. To select the macroscopic reference state, we choose:
\begin{equation}\label{eqn:ref_linear}
U_\ast=U_0^n+\frac{1}{2}(U_0^n-U_0^{n-1})\,,
\end{equation}
or $U_\ast = U_0^1$ at the initial time step. We emphasize we do perform linearization in this approximation and a good linearized approximation does require $F$ close to the selected $M_\ast$. It is beyond the scope of the paper justifying it holds true with our selection in~\eqref{eqn:ref_linear}, and in the numerical example to be demonstrated below, it is clear that when the incoming data is far from the selected local Maxwellian, the approximation breaks down.

The same derivation is used for computing the fluxes at the interface $x=x_{N/2+1/2}=1/2$. Here note that the boundary layer is facing the left, and thus one sets $z = -\frac{x-1/2}{\varepsilon}$ and the velocity also flips the sign:
\begin{equation}\label{eqn:f_boundary_interface}
f(z=0,v) = \frac{F_{N+1}(-v) - M_{\ast}(-v)}{\sqrt{M_{\ast}(-v)}}\,.
\end{equation}
The fluxes for the kinetic region is then given by:
\begin{equation}\label{eqn:F_interface}
F_{N+1,j} = M_\ast(-v_j) + \sqrt{M_\ast(-v_j)}f(z=0,-v_j)\,,\quad\forall v_j>0\,.
\end{equation}

We summarize the algorithm below:

\RestyleAlgo{boxruled}
\begin{algorithm}[H]
\KwData{\begin{itemize}
\item[1] Kinetic boundary condition $F_l(v,t_n)$ for $v>0$ and $F_r(v,t_n)$ for $v<0$;
\item[2] Kinetic solution at $t_n$: $F_{ij}^n$ for $i=N/2+1,\cdots,N$ and all $j$;
\item[3] Fluid solution at $t_n$: $\mathcal{U}_i^n=(\rho_i^n,\rho_i^nu_i^n,E_i^n)^T$ for $i=1\,\cdots N$.\end{itemize}}
\KwResult{Solution at $t_{n+1}$, including: \begin{itemize}
\item[1] Kinetic solution at $t_{n+1}$: $F_{ij}^{n+1}$ for $i=N/2+1,\cdots,N$ and all $j$;
\item[2] Fluid solution at $t_{n+1}$: $\mathcal{U}_i^{n+1}$ for $i=1\,,\cdots N/2$.\end{itemize}}
Step I: prepare boundary fluxes:\\
\qquad Compute~\eqref{eqn:f_layer_nonlinear}, and use~\eqref{eqn:boundary_flux} for $\mathcal{F}^n_{1/2}$;\\
\qquad Compute~\eqref{eqn:f_layer_nonlinear} using the boundary condition in~\eqref{eqn:f_boundary_interface}, and use~\eqref{eqn:boundary_flux} for $\mathcal{F}^n_{N/2+1/2}$,~\eqref{eqn:F_interface} for $F^{n}_{N/2+1}$;\\
Step II: update fluid equation using~\eqref{eqn:Roe method};\\
Step III: update BGK equation using~\eqref{eqn:update_kinetics}.\\
     \caption{Updating the coupled BGK-Euler system from $t_n$ to $t_{n+1}$}\label{alg:euler}
     \end{algorithm}

\subsection{Numerical examples}
We perform a few numerical examples on nonlinear Euler equations with Knudsen layer corrections in this subsection. Throughout the section, for computing the BGK equation, we use $h = 0.001$ and $\Delta t = h/20$. For the numerical velocity range $v\in[-16,16]$ in~\eqref{eqn:vj} we use $M=100$. $\varepsilon$ is set to be $\{\frac{1}{32},\frac{1}{64},\frac{1}{128},\frac{1}{256}\}$. For computing the Euler equations we use coarser grids by setting $h=0.005$ and $\Delta t=0.001$. To measure the error we define
\begin{equation*}
D_\rho=\Vert \rho_\text{BGK}-\rho_\text{Euler}\Vert_{L_2(\text{interior})}\,,
\end{equation*}
and $D_u$ and $D_T$ are defined similarly.

\textbf{Test 4: Pure fluid over the entire domain with boundary layer emerging at $x=0$.} In this example the computed domain is $[0,1]$ with $\varepsilon$ being uniformly small over the entire domain. The initial and boundary conditions for the BGK equation are given as:
\begin{equation*}
  \left\{
    \begin{array}{ll}
      \text{Boundaries: }F_l(v,t)=0,\quad F_r(v,t)=\frac{1}{\sqrt{2\pi}}\exp\Big(\frac{-(v-0.1)^2}{2}\Big) \\
      \text{Initial: }F_0(x,v)=\frac{1}{\sqrt{2\pi}}\exp\Big(\frac{-(v-0.1)^2}{2}\Big)
    \end{array}
  \right. .
\end{equation*}
The initial condition and the right boundary condition are set to be the same and avoid initial and right layers. Correspondingly by using~\eqref{eqn:moments} we obtain the initial condition for the Euler equations:
\begin{equation*}
  \text{Initial: }\rho=1,\quad u=0.1,\quad T=1\,.
\end{equation*}
The boundary condition for the Euler system is computed on-the-fly and cannot be prescribed beforehand.

\begin{figure}[htb]
  \centering
  \includegraphics[width=19cm, height=5cm]{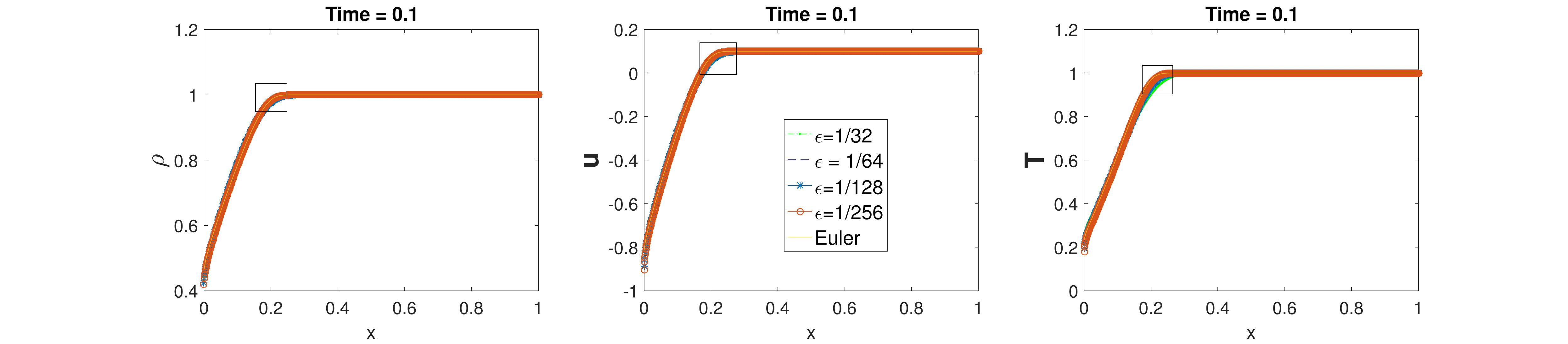}\\
  \caption{Test 4: Solution at Time$=0.1$.}
  \label{fig:Test4}
\end{figure}

\begin{figure}[htb]
  \centering
  \includegraphics[width=19cm, height=5cm]{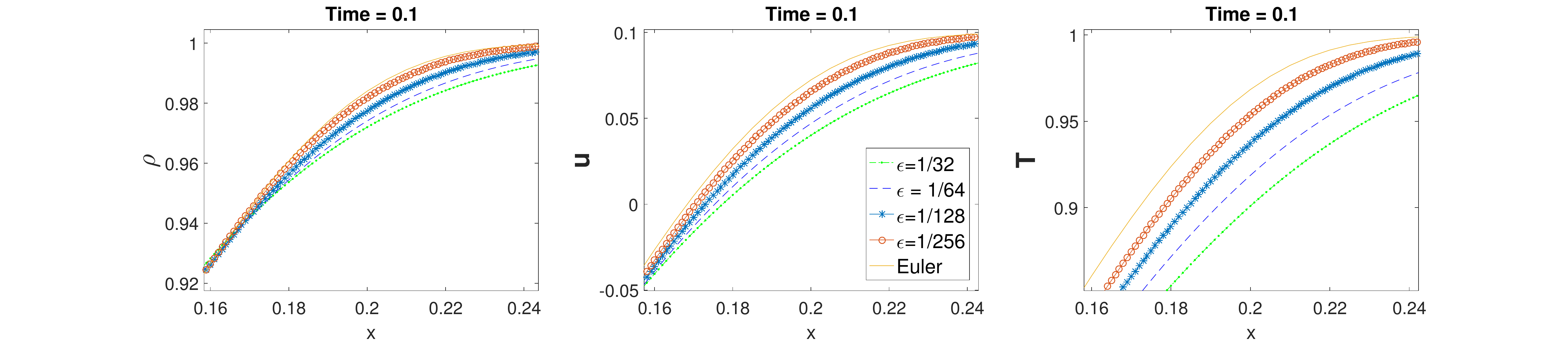}\\
  \caption{Test 4: Zoomed in the box region in Figure~\ref{fig:Test4} .}
  \label{fig:Test4zoom}
\end{figure}

\begin{figure}[htb]
  \centering
  \includegraphics[width=19cm, height=5cm]{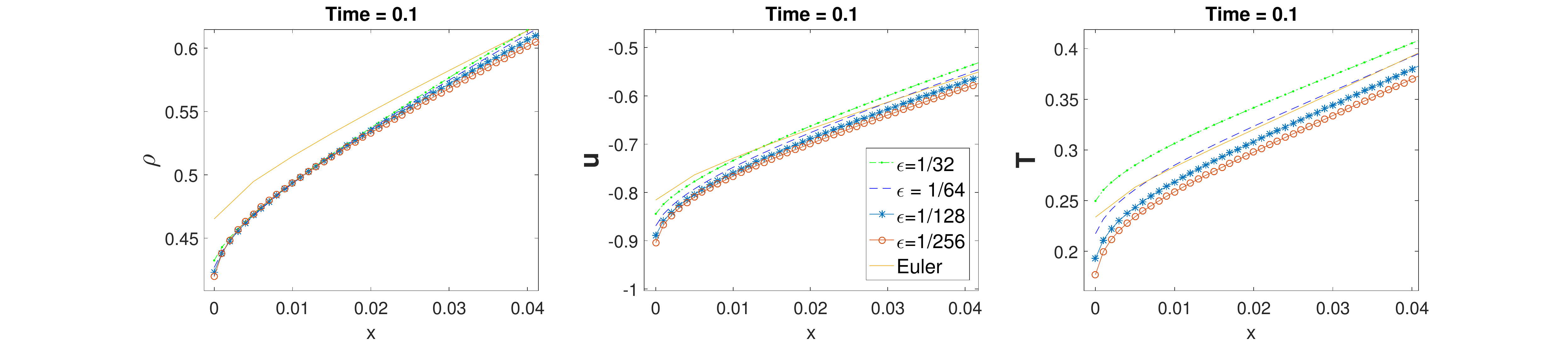}\\
  \caption{Test 4: Boundary layer zoomed at $x=0$.}
  \label{fig:Test4zoom2}
\end{figure}

\begin{figure}[htb]
  \centering
  \includegraphics[width=19cm, height=5cm]{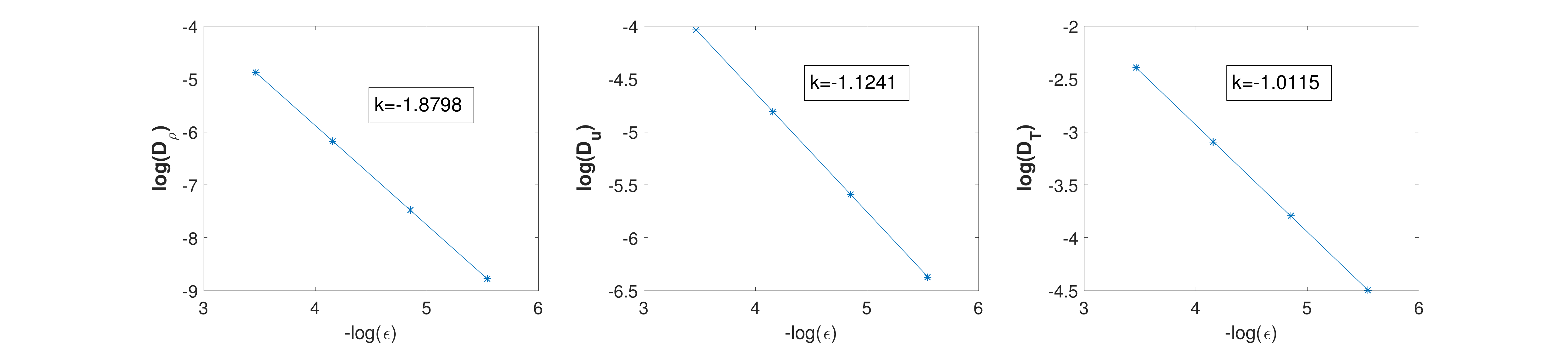}\\
  \caption{Test 4: Error $D_{{\rho}}$, $D_{{u}}$ and $D_{{T}}$ as functions of $\varepsilon$ in log-log scale.}
  \label{fig:Test4error}
\end{figure}

We compute the system up to $\text{Time}=0.1$. In Figure~\ref{fig:Test4} and~\ref{fig:Test4zoom} we plot the solution and the zoom-in of a small region at $\text{Time}=0.1$. Boundary layer emerges at $x=0$, and is shown at Figure~\ref{fig:Test4zoom2}. Even in the nonlinear setting one can see smaller $\epsilon$ leads to closer approximation to the Euler limit. We also show the error decay in a log-log plot in Figure~\ref{fig:Test4error}.

We note that the boundary layer is in fact quite off, if compared with the studies in the linearized setting (\ref{fig:Test1error} for example), especially small $\epsilon$ does not necessarily provides monotonically better approximation to the Euler limit. However, outside the layer zone, the macroscopic quantities are captured rather well. Considering the discrepancy between the boundary data and the interior is relatively big ($\rho\sim 0.4$ at the boundary but $\rho=1$ in the interior), this certainly is an promising evidence that linearization, despite being unsupported of any analytical result, nevertheless provides a good interior solution.

\newpage

\textbf{Test 5: Pure fluid over the entire domain with perturbation on the left boundary.} In this example,  we examine the effect of perturbation by the boundary incoming data to the equilibrium. In particular, we start from a global Maxwellian and add perturbations to the left boundary. The right boundary conditions and the initial condition for the BGK equation are given as:
\begin{equation*}
F_r(v,t)=F_0(x,v)=\frac{1}{\sqrt{2\pi}}\exp\Big(\frac{-(v-0.1)^2}{2}\Big)\,,
\end{equation*}
and we consider two set of boundary conditions on the left boundary, namely
\begin{equation*}
F_l(v,t)^{\text{small}}=\frac{1+5t}{\sqrt{2\pi}}\exp\Big(\frac{-(v-0.1)^2}{2}\Big)\,,
\end{equation*}
and
\begin{equation*}
F_l(v,t)^{\text{large}}=\frac{1+50t}{\sqrt{2\pi}}\exp\Big(\frac{-(v-0.1)^2}{2}\Big)\,,
\end{equation*}
so they deviate from the initial Maxwellian as $t$ grows at a
different rate.

\begin{figure}[htb]
  \centering
  \includegraphics[width=18cm, height=12cm]{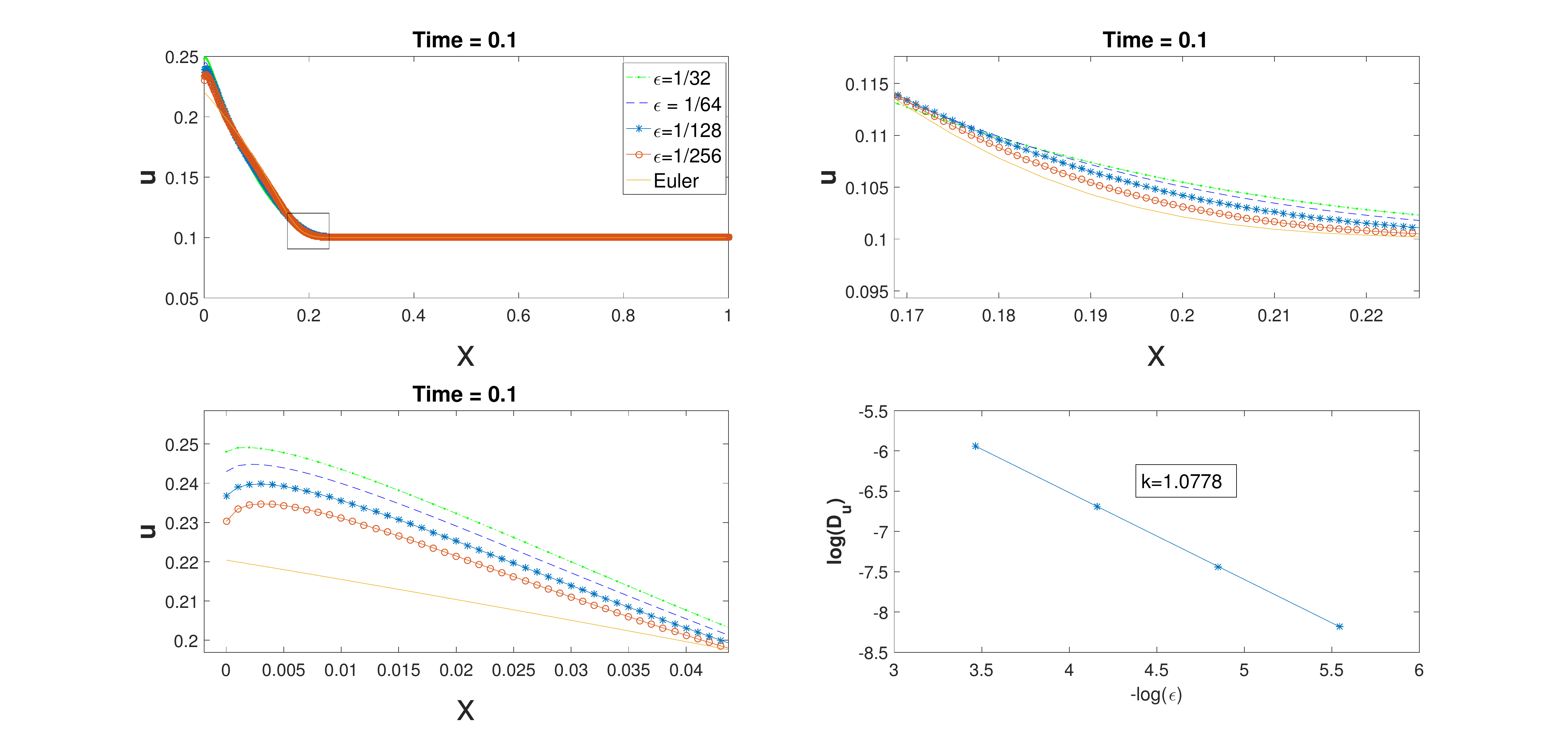}\\
  \caption{Test 5: Solution, zoom-in, error $D_u$ at Time=0.1 for  perturbed boundary data $F_l(v,t)^{\text{small}}$.}
  \label{fig:Test51}
\end{figure}

\begin{figure}[htb]
  \centering
  \includegraphics[width=18cm, height=12cm]{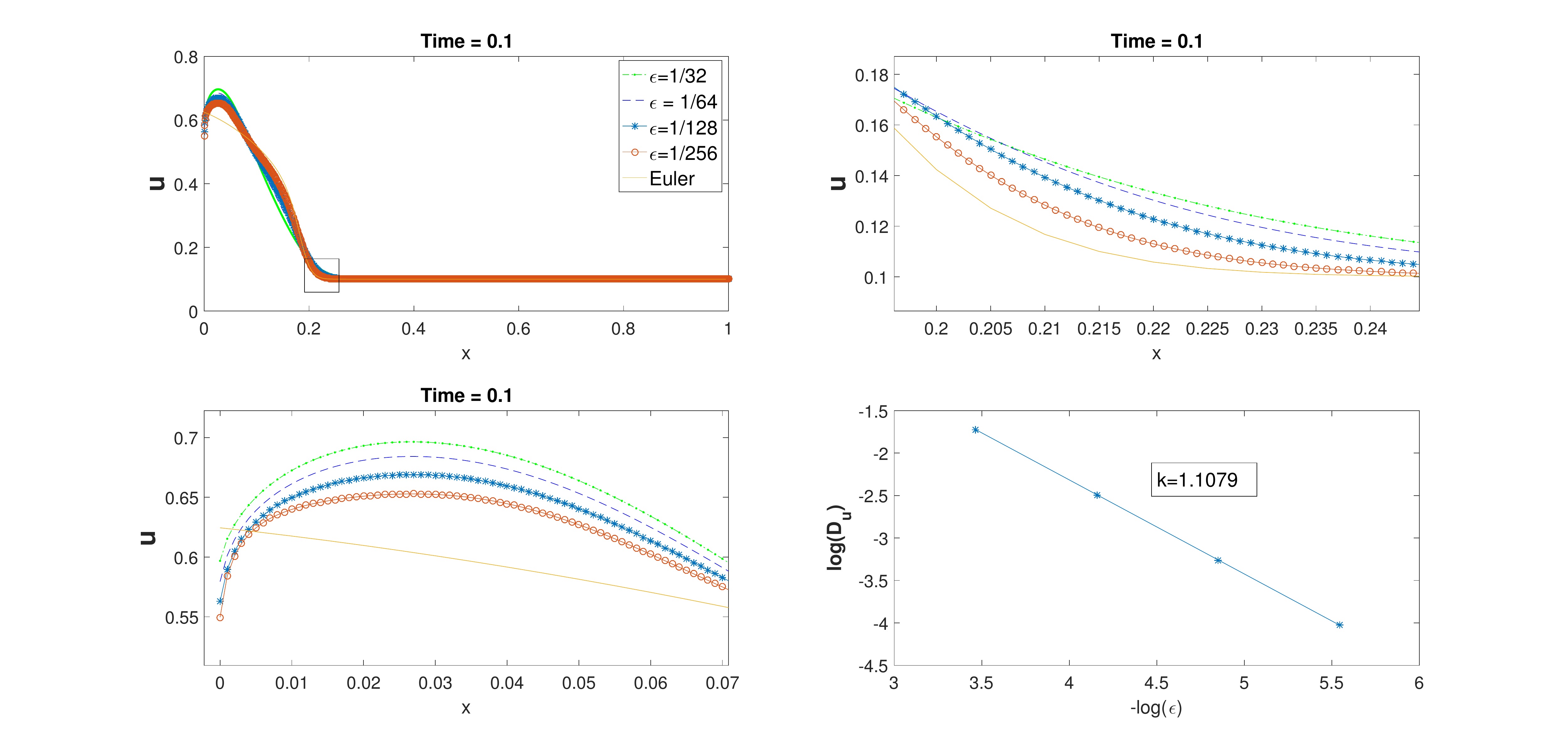}\\
  \caption{Test 5: Solution, zoom-in, error $D_u$ at Time=0.1 for  perturbed boundary data $F_l(v,t)^{\text{large}}$.}
  \label{fig:Test52}
\end{figure}

We compute the system up to $\text{Time}=0.1$. In Figure~\ref{fig:Test51} we plot the solution at $\text{Time}=0.1$, a zoom-in of a small region, zoom-in of the left boundary layer, and the error decay rate of $D_u$ for the $F_l(v,t)^{\text{small}}$ case, and in Figure~\ref{fig:Test52} we plot the counterparts for the $F_l(v,t)^{\text{large}}$ case. We observe, as the previous example, the boundary layer is not well captured, and the is quite far away from the limiting Euler boundary data, but the discrepancy diminishes as the solution propagates into the domain, and the error even still decays with the rate of $\mathcal{O}(\epsilon)$. We also observe that the linearization gives less accurate solution in Figure~\ref{fig:Test52} for the case with a larger perturbation. This is expected since the large perturbation drives the system away from the linearized regime and nonlinear boundary layer effects become more important to capture.

\textbf{Test 6: Coupling the Euler and the BGK equation with interface layer emerging at $x=1/2$.} In this test the fluid limit holds true in the right part of the domain $[0.5, 1]$ in which we use the Euler equations and couple it with the BGK equation computed in $[0,0.5]$. We avoid the boundary layers by setting up compatible initial conditions, and thus only one interal layer emerges at the fluid-kinetic interface at $x=0.5$.

The initial and boundary conditions for the BGK equation are given as:
\begin{equation}\label{}
  \left\{
    \begin{array}{ll}
      \text{Boundaries: } F_l(x,v)=\frac{1}{\sqrt{2\pi}}e^{-(v-0.1)^2/2},\quad F_r(x,v)=\frac{2}{\sqrt{4\pi}}e^{-(v-0.2)^2/4} \\
      \text{Initial: }F_0(x,v)=\left\{
             \begin{array}{ll}
               \frac{1}{\sqrt{2\pi}}e^{-(v-0.1)^2/2}, & \hbox{$x\leq 1/2$;} \\
               \frac{2}{\sqrt{4\pi}}e^{-(v-0.2)^2/4}, & \hbox{$x> 1/2$.}
             \end{array}
           \right.
    \end{array}
  \right.\,
\end{equation}
Correspondingly we use the initial conditions for macroscopic quantities for the Euler equations as:
\begin{equation*}
\rho=\left\{
         \begin{array}{ll}
           1, & \hbox{$x\leq 1/2$;} \\
           2, & \hbox{$x>1/2$.}
         \end{array}
       \right. \quad u=\left\{
                         \begin{array}{ll}
                           0.1, & \hbox{$x\leq 1/2$;} \\
                           0.2, & \hbox{$x>1/2$.}
                         \end{array}
                       \right. \quad T=\left\{
                                         \begin{array}{ll}
                                           1, & \hbox{$x\leq 1/2$;} \\
                                           2, & \hbox{$x>1/2$.}
                                         \end{array}
                                       \right.
\end{equation*}

\begin{figure}[htb]
  \centering
  \includegraphics[width=19cm, height=5cm]{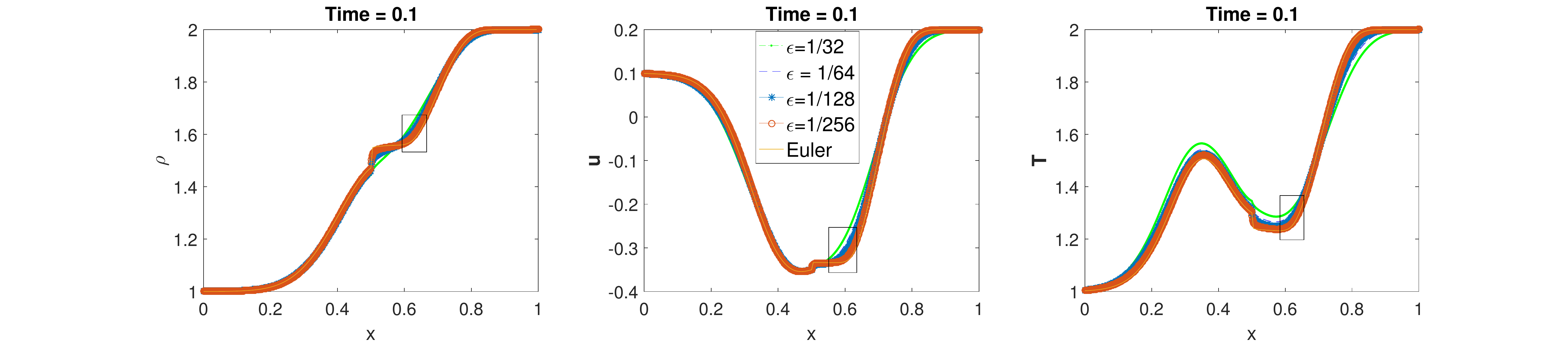}\\
  \caption{Test 6: Solution at Time$=0.1$.}
  \label{fig:Test6}
\end{figure}

\begin{figure}[htb]
  \centering
  \includegraphics[width=19cm, height=5cm]{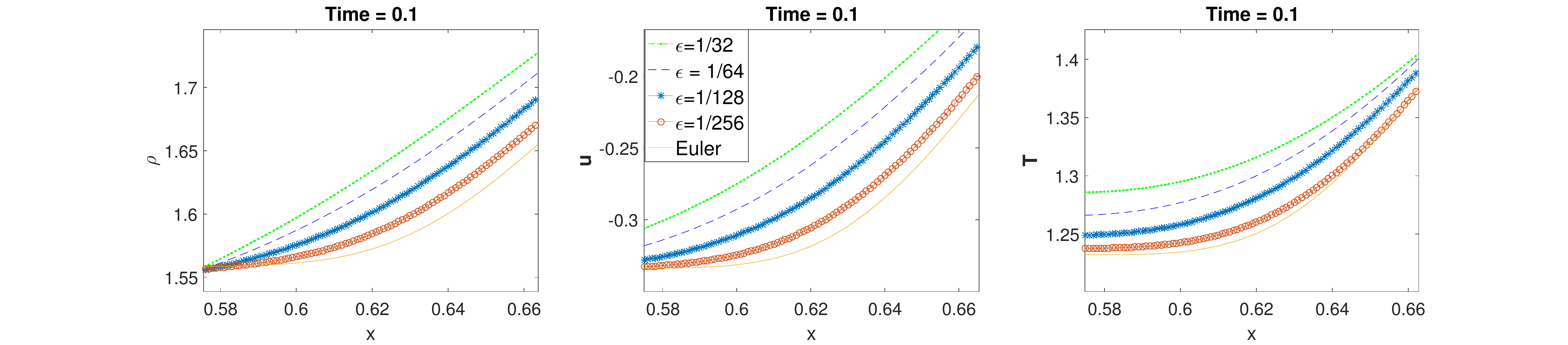}\\
  \caption{Test 6: Zoomed in the box region in Figure~\ref{fig:Test6}.}
  \label{fig:Test6zoom}
\end{figure}

\begin{figure}[htb]
  \centering
  \includegraphics[width=19cm, height=5cm]{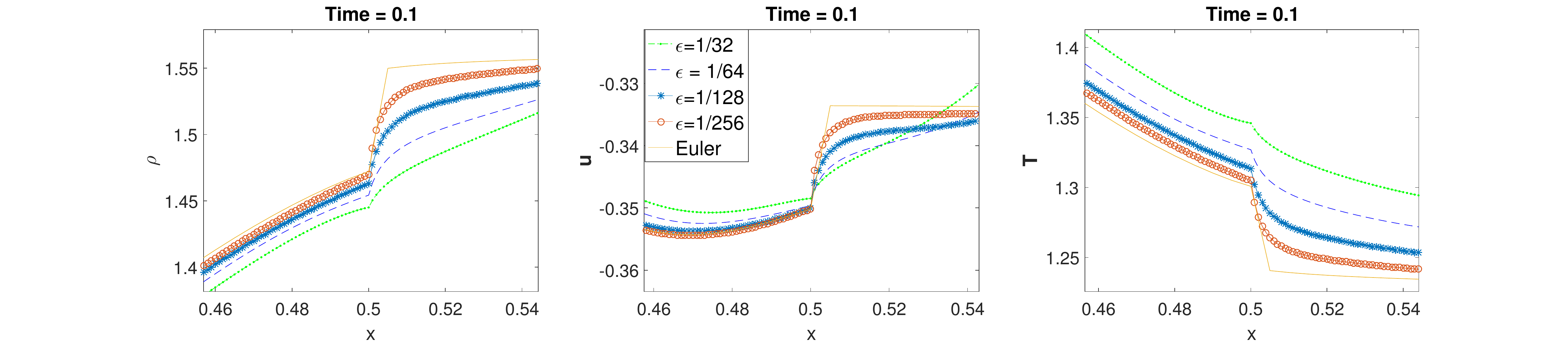}\\
  \caption{Test 6: Interface layer zoomed at $x=0.5$.}
  \label{fig:Test6zoom2}
\end{figure}

\begin{figure}[htb]
  \centering
  \includegraphics[width=19cm, height=5cm]{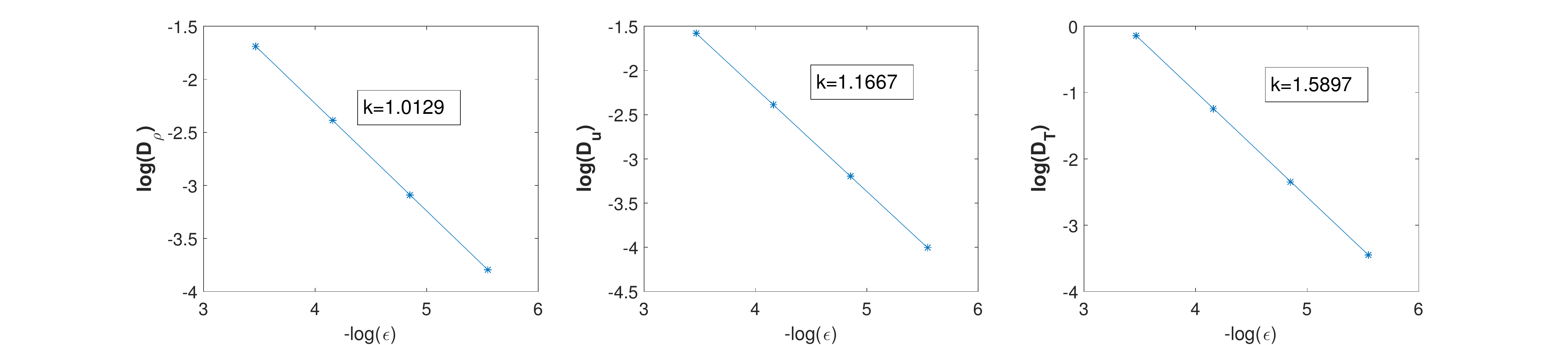}\\
  \caption{Test 6: Error $D_{{\rho}}$, $D_{{u}}$ and $D_{{T}}$ as functions of $\varepsilon$ in log-log scale.}
  \label{fig:Test6error}
\end{figure}

We compute the system up to $\text{Time}=0.1$ using Algorithm~\ref{alg:euler}. In Figure~\ref{fig:Test6}--\ref{fig:Test6error} we plot the solution at the final time, the zoom-in of the interface and the decay of the error terms as $\varepsilon$ decreases. The interface layer is significant as demonstrated in Figure \ref{fig:Test6zoom2}. The proposed method captures the behavior of the interface layer accurately.

\newpage

\bibliographystyle{amsplain}
\bibliography{coupling_BGK}

\appendix

\section{Numerical scheme for half-space kinetic equations}\label{sec:app_A}
In the appendix we present the numerical scheme for half-space kinetic
equations. The numerical method presented in \cite{lls:17_3, lls:15_2}
only considers a fixed reference state $\rho_\ast=1,T_\ast=1/2$, while
the approximation method we developed in this paper involves changing
reference state, thus we need to generalize the numerical methods for
all reference states.
\subsection{Numerical method for boundary layer equation}
In this subsection, we focus on the numerical method for the
half-space problems under the linearized BGK operator. This is a
similar case of the algorithm proposed in the previous work
\cite{lls:17_3, lls:15_2}. Consider the half-space problem, here we
shift the original equation by $u_\ast$,
\begin{equation}\label{}
  \left\{
    \begin{array}{ll}
      (v+u_\ast)\partial_x f+\mathcal{L}f=0, & \hbox{} \\
      f(0,v)=f_0(v), & \hbox{$v+u_\ast >0$,} \\
      f(x,v)\to \theta_{\infty}\in H^0\oplus H^+, & \hbox{$x\to \infty$.}
    \end{array}
  \right.   
\end{equation}
To solve the infinite domain problem, we use a spectral discretization for the $v$-variable. In general the solution may exhibit singularity like jumps at $v=-u_\ast$. Hence we use an even-odd decomposition of the distribution function to avoid the Gibbs phenomena and ensure the accuracy. Here we define the shifted even and odd parts of a function as
\begin{equation}\label{}
  f^E(v)=\frac{f(v)+f(-2u_\ast-v)}{2},\quad f^O(v)=\frac{f(v)-f(-2u_\ast-v)}{2}
\end{equation}
such that $f=f^E+f^O$. Due to the symmetry, it suffices to discretize the function $f^E$ and $f^O$ for $v\in (-u_\ast,\infty)$ and then extend the functions to the whole interval $v\in (-\infty,\infty)$. In other words, we use the half-space general weight Hermite polynomials as basis functions. The construction of the basis functions will be given in the next subsection. Then the even-odd extension is given by
\begin{equation}\label{}
  B_{m}^E(v+u_\ast)=\left\{
               \begin{array}{ll}
                 B_m(v+u_\ast)/\sqrt{2}, & \hbox{$v>-u_\ast$} \\
                 B_m(-v-u_\ast)/\sqrt{2}, & \hbox{$v<-u_\ast$}
               \end{array}
             \right.
\quad B_m^O(v+u_\ast)=\left\{
                                     \begin{array}{ll}
                                       B_m(v+u_\ast)/\sqrt{2}, & \hbox{$v>-u_\ast$} \\
                                       -B_m(-v-u_\ast)/\sqrt{2}, & \hbox{$v<-u_\ast$}
                                     \end{array}
                                   \right. 
\end{equation}
where $B_m(v+u_\ast)$ are the general weight Hermite polynomials on $(-u_\ast,\infty)$ satisfying
\[\int_{-u_\ast}^\infty B_n(v+u_\ast)B_m(v+u_\ast)e^{-(v+u_\ast)^2/2T}dv=\int_0^\infty B_n(v)B_m(v)e^{-v^2/2T}dv=\delta_{nm}.\]
Finally the basis functions $P_n$ are obtained by multiplying these functions by the square root of the Maxwellian:
\begin{eqnarray}
  P_{2n-1}(v+u_\ast) &=& B_{n-1}^O(v+u_\ast) e^{-(v+u_\ast)^2/2T}, \\
  P_{2n}(v+u_\ast) &=& B_{n-1}^E(v+u_\ast) e^{-(v+u_\ast)^2/2T}.
\end{eqnarray}
For the stability of the numerical method, we first solve a damped version of and then recover the solution to the original equation. Note that after the shift the basis function of the null space of $\mathcal{L}_\ast$ is given by
\begin{equation}\label{}
  \left\{
    \begin{array}{ll}
      \chi_0(v)=\frac{1}{\sqrt{6}}(\frac{v^2}{T_\ast}-3)\sqrt{M_{[0,T_\ast]}}, &  \\
      \chi_+(v)=\frac{1}{\sqrt{6}}(\sqrt{\frac{3}{T_\ast}}v+\frac{v^2}{T_\ast})\sqrt{M_{[0,T_\ast]}}, &  \\
      \chi_-(v)=\frac{1}{\sqrt{6}}(\sqrt{\frac{3}{T_\ast}}v-\frac{v^2}{T_\ast})\sqrt{M_{[0,T_\ast]}}. &
    \end{array}
  \right. 
\end{equation}
The damped equation is given by
\begin{equation}\label{}
  \left\{
    \begin{array}{ll}
      (v+u_\ast)\partial_x f+\mathcal{L}_df=0, & \hbox{} \\
      f(0,v)=f_0(v), & \hbox{$v+u_\ast >0$,}

    \end{array}
  \right.
\end{equation}
where
\begin{multline}\label{}
 \mathcal{L}_df=\mathcal{L}f+\sum_{k=1}^{v_+}\alpha (v+u_\ast)\chi_+\langle (v+u_\ast)\chi_+,f\rangle+\sum_{k=1}^{v_0}\alpha (v+u_\ast)\chi_0\langle (v+u_\ast)\chi_0,f\rangle \\
 +\sum_{k=1}^{v_-}\alpha (v+u_\ast)\chi_-\langle (v+u_\ast)\chi_-,f\rangle +\sum_{k=1}^{v_0}\alpha (v+u_\ast)\mathcal{L}^{-1}((v+u_\ast)\chi_0)\langle (v+u_\ast)\mathcal{L}^{-1}((v+u_\ast)\chi_0)\rangle .
\end{multline}
The well-posedness of this equation is proved in Proposition 3.2 \cite{lls:17_3}, which verifies the inf-sup condition of the variational formulation. We approximate the even and odd parts of the distribution functions by
\begin{equation}\label{}
  f^E(x,v)=\sum_{n=1}^N a_n^E(x)P_n^E(v),\quad f^O(x,v)=\sum_{n=1}^{N+1}a_n^OP_n^O(v).
\end{equation}
Substituting the approximation into and applying Galerkin method, we obtain the equation for the coefficients which reads
\begin{equation}\label{}
  A\partial_x \vec{a}=B\vec{a},
\end{equation}
where
\begin{equation}\label{}
  A_{ij}=\langle vP_i,P_j\rangle,\quad B_{ij}^E=\langle \mathcal{L}_dP_i,P_j\rangle.
\end{equation}
After diagonalizing the equation into a generalized eigenvalue problem, we obtain a system of $2N+1$ ODE reads
\begin{equation}\label{}
  \partial_x \vec{b}=V\vec{b},
\end{equation}
with $A^{-1}B=XVX^{-1}$ and $b=Xa$. $V$ is a diagonal matrix. The solution of the ODE tell us that we need $2N+1$ boundary conditions to determine $\vec{b}$. The boundary conditions are of two kinds. The first is given by the Dirichlet boundary condition. Note that the boundary condition only provides data at $v>-u_\ast$, we only get $N$ conditions for $\vec{b}$. The remaining conditions come from the requirement that $f\to \theta_{\infty}\in H^0\oplus H^+$, this means $\vec{a}$ can not be exponential increasing. Hence positive eigenvalue corresponds to 0 coefficient of $\vec{b}$. It is proved in Proposition 4.6 \cite{lls:17_3} that there are exactly $N$ positive eigenvalues and 1 zero eigenvalue of the generalized eigenvalue. We obtain enough conditions to determine $\vec{a}=X\vec{b}$.
\\Once we obtain the solution of the damped equation, we can explicitly construct solutions to the undamped equation as stated in Theorem \ref{thm:boundary_compute}. Specifically, let $g_+$ be the solution to (47) with boundary conditions given by $\chi_+$:
\[g_+|_{x=0}=\chi_+,\quad v+u_\ast>0.\]
Similarly, denote $g_0$ as the solution to (47) where the incoming boundary data is given by $\chi_0$. Let $C$ be the block matrix defined by
\begin{equation}\label{}
  C=\left(
      \begin{array}{cc}
        C_{++} & C_{+0} \\
        C_{0+} & C_{00} \\
      \end{array}
    \right),
\end{equation}
where
\begin{eqnarray*}
  C_{++}=\langle (v+u_\ast)\chi_+,g_+\rangle|_{x=0} ,\quad C_{+0}=\langle (v+u_\ast)\chi_{+},g_0\rangle|_{x=0} \\
   C_{0+}=\langle (v+u_\ast)\chi_0,g_+\rangle|_{x=0} ,\quad C_{00}=\langle (v+u_\ast)\chi_{0},g_0\rangle|_{x=0}.
\end{eqnarray*}
Define the coefficient vector $\eta=(\eta_+,\eta_1)^T$ such that
\begin{equation}\label{}
  C\eta=U_f,
\end{equation}
where $U_f=(u_+,u_0)$ with $u_+=\langle (v+u_\ast)\chi_+,f\rangle_{x=0}$ and $u_0=\langle (v+u_\ast)\chi_0,f\rangle_{x=0}$.
\\
In fact, $C$ is invertible and hence (66) is uniquely solvable, moreover,
\begin{equation}\label{}
  f_{\phi}=f-\sum_{k=1}^{v_+}\eta_+(g_+-\chi_+)-\sum_{k=1}^{v_0}\eta_0(g_0-\chi_0)
\end{equation}
is the unique solution to the half-space equation
\begin{eqnarray}
  (v+u_\ast)\partial_x f+\mathcal{L}f &=& 0 ,\\
  f|_{x=0} &=& \phi(v)\quad v+u_\ast>0.
\end{eqnarray}
Moreover, the end state $f_{\phi,\infty}$ is given by
\[f_{\phi,\infty}=\sum_{k=1}^{v_+}\eta_+\chi_++\sum_{k=1}^{v_0}\eta_0\chi_0.\]
In sum, we use the algorithm described above to obtain the solution of the half-space kinetic equation and then use the algorithm described in previous section to deal with the coupling problems.

\subsection{General half-space Hermite polynomials}
As described in previous subsection, to generalized the algorithm for general reference state, we can see the key is to generalized the half-space Hermite polynomials. Then the Galerkin method remains the same.

The basis functions is constructed using the half-space Hermite polynomials, which are orthogonal polynomials defined on the positive half $v$-axis with the weight function $e^{-v^2/2T}$: $\{B_n(v),v>0\}$ such that each $B_n(v)$ is a polynomial of order $n$ and
\begin{equation}\label{}
  \int_0^\infty B_m(v)B_n(v)e^{-v^2/2T}dv=\delta_{nm},
\end{equation}
The orthogonal polynomials can be constructed using three term recursion formula, for the derivation one can see the details in appendix.

The basis function we need are either odd or even with respect to $v=\-u$, thus we shift $B_n$ by $-u$ and make even-odd extension
\begin{eqnarray}
  B_n^E(v) &=& \left\{
                 \begin{array}{ll}
                   B_n(v+u)/\sqrt{2}, & \hbox{$v>-u$;} \\
                   B_n(-v-u)/\sqrt{2}, & \hbox{$v<-u$.}
                 \end{array}
               \right.
 \\
  B_n^O(v) &=& \left\{
                 \begin{array}{ll}
                   B_n(v+u)/\sqrt{2}, & \hbox{$v>-u$;} \\
                   -B_n(-v-u)\sqrt{2}, & \hbox{$v<-u$.}
                 \end{array}
               \right.
\end{eqnarray}
The $(2N+1)\times (2N+1)$ matrices are given by
\[A_{ij}=\int_{\mathbb{R}}(v+u)P_iP_jdv,\quad B_{ij}=-\int_{\mathbb{R}}P_i\mathcal{L}_dP_jdv.\]
$A$ can be obtained by the recurrence relation. For matrix $B$, recall that

\begin{multline}\label{}
 \mathcal{L}_df=\mathcal{L}f+\sum_{k=1}^{v_+}\alpha (v+u)\chi_+\langle (v+u)\chi_+,f\rangle+\sum_{k=1}^{v_0}\alpha (v+u)\chi_0\langle (v+u)\chi_0,f\rangle \\
 +\sum_{k=1}^{v_-}\alpha (v+u)\chi_-\langle (v+u)\chi_-,f\rangle +\sum_{k=1}^{v_0}\alpha (v+u)\mathcal{L}^{-1}((v+u)\chi_0)\langle (v+u)\mathcal{L}^{-1}((v+u)\chi_0)\rangle .
\end{multline}
All the integrals involved in calculating $B$ can be obtained by using the Gaussian quadrature. To see this, we take $\chi_0$ as an example, all other integrals can be treated in the same way. We firstly split the integral into two parts
\[\int_{-\infty}^\infty P_{2i}\chi_0dv=\int_{-u}^\infty P_{2i}\chi_0dv+\int_{-\infty}^{-u}P_{2i}\chi_0dv.\]
Note that $P_{2i}$, on each side of $-u$, is $i$-th order polynomial product with $\exp(-\frac{(v+u)^2}{4T})$ while $\chi_0$ is a quadratic function multiplied with a different weight function $\exp(-v^2/4T)$. The two Gaussians that entered at different locations could be combined, and the numerical integral is exact(for polynomials up to $N$-th order) once the correct Gaussian quadratures is adopted:
\begin{multline}\label{}
  \int_{-u}^\infty  P_n(v)\chi_0(v)dv=\frac{1}{\sqrt{2}}\int_{-u}^\infty B_n(v+u)\chi_0e^{-\frac{(v+u)^2+v^2}{4T}}dv\\
=\frac{1}{\sqrt{2}}\int_0^\infty B_n(v)\chi_0(v-u)e^{-\frac{2v^2-2vu+u^2/2+u^2/2}{4T}}dv=\frac{1}{\sqrt{2}}e^{-u^2/8T}\int_0^\infty B_n(v-u)\chi_0(v-u)e^{\frac{-(v-u/2)^2}{2T}}.
\end{multline}

Similarly, we have the integration from the negative part,
\[\int_{-\infty}^{-u}P_n\chi_0dv=\frac{1}{\sqrt{2}}\int_{-\infty}^{-u}B_n(v+u)\chi_0e^{-\frac{(v+u)^2+v^2}{4T}}dv=\frac{1}{\sqrt{2}}e^{-u^2/8T}\int_0^\infty B_n(-v-u)\chi_0(-v)e^{\frac{-(v+u/2)^2}{2T}}dv.\]
Thus these integrals can be obtained by using the Gaussian quadrature of the weight $e^{-(v-u/2)^2/2T}$ and $e^{-(v+u/2)^2/2T}$ respectively.

\subsection{Derivation of recurrence relation}
Here we derive the half-space orthogonal polynomials with weight $\exp((v-u)^2/2T)$. The zeroth order is
\[B_0=\frac{1}{\sqrt{m_0}},\quad m_0=\frac{\sqrt{\pi}}{2}\Big(1+\textrm{erf}(\frac{u}{2T})\Big).\]
Set the recurrence relation for the higher order polynomials as
\begin{equation}\label{}
  \sqrt{\beta_{n+1}}B_{n+1}=(v-\alpha_n)B_n-\sqrt{\beta_n}B_{n-1},
\end{equation}
we aim to derive the formula for $\beta_n$ and $\alpha_n$. Actually
\begin{equation}\label{}
  \left\{
    \begin{array}{ll}
      \beta_{n+1}=2Tn+T-\beta_n+u\alpha_n-\alpha_n^2 \\
     \alpha_{n+1}=\frac{T}{\beta_{n+1}}\sum_{k=0}^n \alpha_k-\alpha_n+u
    \end{array}
  \right.
\end{equation}
with $\alpha_0=m_1/m_0$, $\beta_0=0$ where $m_i$, i=0,1 are moments of the Gaussian:
\begin{equation}\label{}
  m_i=\int_0^\infty v^i e^{-(v-u)^2/2T}dv.
\end{equation}
Now we start the derivation. From the recurrence relation we get
\[\alpha_n=\int_0^\infty vB_n^2 e^{\frac{-(v-u)^2}{2T}}dv,\quad \sqrt{\beta_{n+1}}=\int_0^\infty vB_nB_{n+1}e^{\frac{-(v-u)^2}{2T}}dv.\]
By the Christoffel-Darboux identity
\begin{equation}\label{}
\sum_{k=0}^n B_k^2=\sqrt{\beta_{n+1}}(B_{n+1}'B_n-B_{n+1}B_n').
\end{equation}
By integrating this identity with the weight we get
\[n+1=\sqrt{\beta_{n+1}}\int_0^\infty B_{n+1}'B_n e^{\frac{-(v-u)^2}{2T}}dv.\]
Note that
\[\sqrt{\beta_{n+1}}B_{n+1}'=B_n+(v-\alpha_n)B_n'-\sqrt{\beta_n}B_{n-1}'.\]
We have
\[n+1=1+\int_0^\infty (v-\alpha_n)B_n' B_n e^{\frac{-(v-u)^2}{2T}}dv=1+\int_0^\infty vB_n' B_n e^{\frac{-(v-u)^2}{2T}}\]
\[n=-\int_0^\infty \frac{1}{2}(B_n)^2 \Big[ve^{\frac{-(v-u)^2}{2T}}\Big]'dv=-\int_0^\infty \frac{1}{2}(B_n)^2 \Big[e^{\frac{-(v-u)^2}{2T}}-\frac{v-u}{T}ve^{\frac{-(v-u)^2}{2T}}\Big]dv\]

\[=-\frac{1}{2}+\int_0^\infty \frac{1}{2T}v^2B_n^2e^{\frac{-(v-u)^2}{2T}}dv-\frac{1}{2T}u\alpha_n.\]
Note that
\[vB_n=\sqrt{\beta_{n+1}}B_{n+1}+\alpha_n B_n+\sqrt{\beta_n}B_{n-1},\]
we get
\[v^2B_n^2=\sqrt{\beta_{n+1}}vB_{n+1}B_n+v\alpha_n B_n^2+v\sqrt{\beta_n}B_nB_{n-1}.\]
Therefore,
\[n=-\frac{1}{2}+\frac{1}{2T}(\beta_{n+1}+\alpha_n^2+\beta_n-u\alpha_n) \]
\[\beta_{n+1}=2Tn+T-\beta_n+u\alpha_n-\alpha_n^2.\]

Next we multiply the identity with $v$ and then integrate to obtain
\[\sum_{k=0}^n \alpha_k=\sqrt{\beta_{n+1}}\int_0^\infty v B_{n+1}'B_n e^{\frac{-(v-u)^2}{2T}}dv\]
\[=\sqrt{\beta_{n+1}}\Big(\int_0^\infty \frac{v^2}{T}B_{n+1}B_n e^{\frac{-(v-u)^2}{2T}}-\int_0^\infty \frac{vu}{T}B_{n+1}B_n e^{\frac{-(v-u)^2}{2T}}dv \Big)=\frac{\beta_{n+1}}{T}(\alpha_n+\alpha_{n+1}-u)\]
\[\alpha_{n+1}=\frac{T}{\beta_{n+1}}\sum_{k=0}^n \alpha_k-\alpha_n+u.\]

\end{document}